\documentclass[final,3p]{elsarticle}
\usepackage{amssymb,amsmath,amsthm,epsfig,verbatim,url}
\usepackage{epstopdf}
\biboptions{sort&compress}

\usepackage[notref,notcite]{showkeys}
\usepackage{mathrsfs}
\usepackage{pgfplots}
\usepackage{ulem}
\usepackage{caption}
\pgfplotsset{compat=1.7}

\newtheorem{thm}{Theorem}[section]
\newtheorem{lem}[thm]{Lemma}
\newtheorem{rem}[thm]{Remark}

\newcommand*{\vv}[1]{\vec{\mkern0mu#1}}

\newcommand{\norm}[1]{\Vert#1\Vert}

\newcommand{\bR}{{\mathbb R}}
\newcommand{\bN}{{\mathbb N}}

\newcommand{\tD}{\mathbb{D}}
\newcommand{\mX}{\mathscr{X}}
\newcommand{\vecz}{\vec{\mathfrak{x}}}

\newcommand{\vol}{\operatorname{vol}}

\newcommand{\dL}{{\rm d}\mathscr{L}}
\newcommand{\dH}{{\rm d}\mathscr{H}}
\newcommand{\mL}{\mathscr{L}}
\newcommand{\mH}{\mathscr{H}}
\newcommand{\rd}{\;{\rm d}}
\newcommand{\id}{{\rm id}}
\newcommand{\dd}[1]{\frac{\rm d}{{\rm d}#1}}
\newcommand{\ddt}{\dd{t}}
\newcommand{\nn}{\nonumber}
\newcommand{\ttau}{\Delta t}

\newcommand{\ESP}{{\rm Eulerian}}
\newcommand{\ALE}{{\rm ALE}}

\newcommand{\cira}{\mbox{$c\!\!\!\!\:/$}}

\newcommand{\mat}[1]{\uuline{#1}}
\newcommand{\Id}{{I\!d}}

\begin{document}
\begin{frontmatter}
\title{
Structure-preserving discretizations of two-phase Navier--Stokes flow using fitted and unfitted approaches
}

\author[1]{Harald Garcke}
\address[1]{Fakult{\"a}t f{\"u}r Mathematik, Universit{\"a}t Regensburg, 
93040 Regensburg, Germany}
\ead{harald.garcke@ur.de}
\author[2]{Robert N\"urnberg}
\address[2]{Dipartimento di Mathematica, Universit\`a di Trento,
38123 Trento, Italy}
\ead{robert.nurnberg@unitn.it}
\author[1]{Quan Zhao}
\ead{quan.zhao@ur.de}

\begin{abstract}
We consider the numerical approximation of a sharp-interface model for 
two-phase flow, which is given by the incompressible Navier--Stokes equations 
in the bulk domain together with the classical interface conditions on the 
interface. We propose structure-preserving finite element methods for the 
model, meaning in particular that volume preservation and energy decay are 
satisfied on the discrete level. 
For the evolving fluid interface, we employ parametric finite element 
approximations that introduce an implicit tangential velocity to improve the 
quality of the interface mesh. For the two-phase Navier--Stokes equations, 
we consider two different approaches: an unfitted and a fitted finite element 
method, respectively. In the unfitted approach, the constructed method is 
based on an Eulerian weak formulation, while in the fitted approach a novel 
arbitrary Lagrangian--Eulerian (ALE) weak formulation is introduced. 
Using suitable discretizations of these two formulations, we introduce 
two finite element methods and prove their structure-preserving properties.
Numerical results are presented to show the accuracy and efficiency of the 
introduced methods. 
\end{abstract} 

\end{frontmatter}

\renewcommand{\thefootnote}{\arabic{footnote}}
\setlength{\parindent}{2em}


\setcounter{equation}{0}
\section{Introduction} \label{sec:intro}

The problem of two-phase flows has attracted a lot of attention in recent 
decades, not only because it involves many interesting phenomena in nature, but also 
due to its important applications in various fields, such as ink-jet printing, 
coating and microfluidics in industrial engineering and scientific experiments.
Therefore, developing accurate and robust numerical methods for these flows is 
necessary and meaningful. 

According to the treatment of the moving interface between the two phases, 
numerical approximations for two-phase flows can be classified into two main 
categories. The first category is based on interface capturing methods, 
where the interface is determined implicitly by an auxiliary scalar function 
defined on a fixed  domain. These include the volume of fluid method 
\cite{Hirt1981volume,Renardy02,Popinet09}, the level set method 
\cite{Sussman94level,Sethian99level,Gross07extended,osher02level,olsson07}, 
and the diffuse-interface method 
\cite{Anderson1998,Feng2006fully,Styles2008finite,Abels2012,Grun2014two, Aland2012}. 
The second category comprises the so-called front-tracking methods. Here the
interface is explicitly tracked by a collection of markers or by a lower 
dimensional moving mesh, see, e.g., 
\cite{Hughes81,Tryg01,Bansch01,Perot2003moving,Ganesan06,BGN2013eliminating,%
BGN15stable,Agnese20}. 
In general, interface-capturing methods can automatically handle possibly 
complex topological changes in the evolution of the interface. 
In front-tracking methods, on the other hand, topological changes need to be
performed heuristically. In addition, front-tracking methods may struggle
with the preservation of the mesh quality of the moving interface, especially 
during strong deformations of the interface. Nevertheless, front-tracking 
methods offer very accurate and efficient approximations of the interface and 
its geometry. For example, the curvature of the interface, denoted by 
$\varkappa$, can be accurately computed with the help of the identity 
\cite{Dziuk90,Deckelnick2005}  
\begin{equation}
\varkappa \,\vec\nu = \Delta_s\vec\id,
\label{eq:CF}
\end{equation}
where $\vec\id$ is the identity function, $\vec\nu$ is the unit normal to the 
interface, $\Delta_s=\nabla_s\cdot\nabla_s$ is the Laplace--Beltrami operator 
on the interface with $\nabla_s$ being the surface gradient. The identity 
\eqref{eq:CF} was initially used for the computation of mean curvature flow in 
\cite{Dziuk90}, and then has been generalized to approximate the surface 
tension force in the context of two-phase flow, e.g., 
\cite{Bansch01,Ganesan07,BGN2013eliminating,BGN15stable,Agnese16,Agnese20,%
Zhao2020energy,Zhao2021,Duan2022energy}.

 Typically numerical methods for two-phase flow are such that the choice of how to approximate the moving interface is completely independent from the choice of how to discretize the flow in the bulk.
But for front-tracking methods a design choice has to be made.
In fact, the latter methods can be divided into fitted and
unfitted, distinguishing between two possible
relationships between the approximation of the moving interface and the
bulk mesh. In the unfitted approach, the bulk mesh and the interface mesh are totally 
independent, thus allowing the interface to cut through the elements of the 
bulk mesh. One of the prominent examples for unfitted approximations is the 
immersed boundary methods \cite{Peskin2002,Leveque97,Li2013three}. In the 
finite element framework, usually an enrichment of the elements that are cut 
by the interface is necessary in order to accurately capture jumps of physical
quantities across the interface, see XFEM 
\cite{Gross07extended,Ausas2012new,BGN2013eliminating,BGN15stable} and 
cutFEM \cite{Frachon2019cut,Claus2019cutfem}. 
In the fitted mesh approach, the discrete interface is made up of faces of 
elements from the bulk mesh, and thus the bulk mesh needs to deform 
appropriately in time in order to match the evolving interface. 
Employing Eulerian schemes on this moving bulk mesh implies that the obtained 
solutions need to be frequently interpolated on the new mesh. In purely
Lagrangian schemes, on the other hand, the bulk mesh needs to deform 
according to the fluid velocity, and this often leads to large distortions
of the mesh. A possible way to overcome these drawbacks is the so-called 
arbitrary Lagrangian--Eulerian (ALE) approach, e.g., 
\cite{Hughes81,Ganesan06,Gerbeau2006,Agnese20,Duan2022energy,Anjos20143d}, 
where the equations are formulated in a moving frame of reference, and the 
corresponding reference velocity is independent from the fluid velocity, and in
this sense it is ``arbitrary''. 

{From} the numerical analysis point of view, it is of great interest to 
consider numerical approximations that can preserve the energy-diminishing and 
volume-preserving structure of the considered flow. Based on \eqref{eq:CF}, 
B\"ansch proposed a space-time finite element discretization for free 
capillary flows \cite{Bansch01}, which yields a nonlinear scheme that satisfies a stability bound. Barrett, Garcke and N\"urnberg introduced a novel weak
formulation \cite{Barrett20}, which is referred to as the BGN formulation 
from now on. In this formulation, the interface is advected in the normal 
direction according to the normal part of the fluid velocity, thus allowing 
tangential degrees of freedom to improve the mesh quality. The formulation was 
employed for two-phase Stokes flow \cite{BGN2013eliminating} and for two-phase 
Navier--Stokes flow \cite{BGN15stable} with unfitted finite element 
approximations, which leads to an unconditional stability bound. Moreover, 
applications to the two-phase Navier--Stokes flow with moving fitted finite 
element methods were also investigated in \cite{Agnese20} in both Eulerian and 
ALE approaches, but in general no stability bound was available. Recently, 
by assuming that the ALE frame velocity satisfies the divergence-free 
condition, Duan, Li and Yang proposed an ALE weak formulation, which leads to 
an energy-diminishing scheme on the discrete level \cite{Duan2022energy}. 
In addition to these energy stable approximations, numerical approximations 
that maintain the volume preservation of the two phases are also desirable. 
By enriching the pressure space with extra degrees of freedom 
\cite{BGN2013eliminating,BGN15stable,Agnese16}, the finite element 
approximations derived from the BGN formulation can achieve the volume 
preservation on the semidiscrete level. However, on the fully discrete level
an exact volume preservation in general does not hold. 

Recently, based on the BGN formulation \cite{Barrett20} and the idea in 
\cite{Jiang21}, Bao and Zhao \cite{BZ21SPFEM} proposed a numerical method for 
surface diffusion flow, which enables exact volume conservation for the fully 
discrete solutions with the help of time-integrated discrete normals. 
In this paper we would like to incorporate this novel idea into the numerical
approximation of two-phase flows.
For completeness we note that in \cite{Li2013three}, a volume-correction 
method was proposed, which works by relocating the interface points along the 
normal direction of the interface in a suitable way. However, to our 
knowledge, no existing methods can exactly preserve the volume of the two 
phases intrinsically on the fully discrete level.

The main aim of this work is to develop structure-preserving discretizations 
for two-phase Navier--Stokes flow in both the unfitted and fitted mesh 
approaches. In the unfitted approach, we combine the ideas in 
\cite{BGN15stable} and \cite{BZ21SPFEM} to obtain a fully discrete 
approximation that satisfies an unconditional stability estimate and an
exact volume preservation property. In the fitted mesh approach, we employ a
novel ALE moving mesh formulation together with the idea in \cite{BZ21SPFEM}. 
The introduced scheme is based on suitable discretizations of a novel 
nonconservative ALE weak formulation. We argue that in terms of the treatment 
of the inertia term, our ALE method is similar to the work in 
\cite{Duan2022energy}, but here we allow for a more flexible mesh velocity of 
the ALE frame, which helps to guarantee the quality of the bulk mesh.  

The rest of the paper is organized as follows. We start in 
Section~\ref{sec:mathF} by reviewing the strong formulation of two-phase 
Navier--Stokes flow. Next, in Section~\ref{sec:unfem} we are focused on the 
unfitted mesh approach with an Eulerian weak formulation. This leads to a 
``weakly'' nonlinear discretized scheme that enjoys volume conservation and 
unconditional energy stability. In Section~\ref{sec:fitfem}, we move to the 
fitted mesh approach. A structure-preserving method is introduced based on 
suitable discretizations of a novel ALE weak formulation. Subsequently, we 
present several benchmark tests for the introduced schemes in 
Section~\ref{sec:num}. Finally, we draw some conclusions in 
Section~\ref{sec:con}. 

\setcounter{equation}{0}
\section{The strong formulation}\label{sec:mathF}

\begin{figure}[!htp]
\centering
\includegraphics[width=0.45\textwidth]{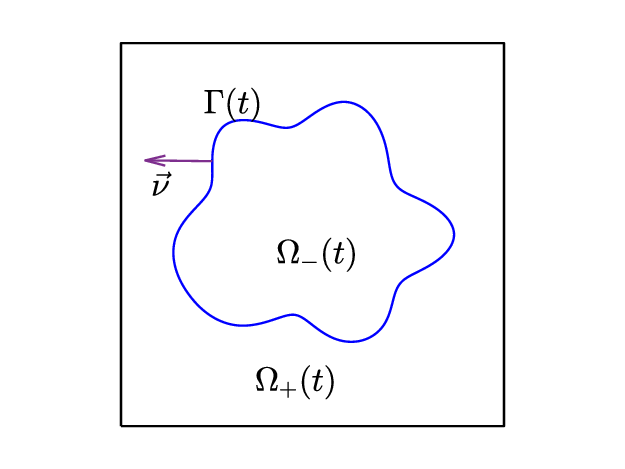}
\caption{An illustration of two-phase flow in a bounded domain $\Omega=\Omega_-(t)\cup\Omega_+(t)\cup\Gamma(t)$ in the case $d=2$.}
\label{fig:tpf}
\end{figure}

As shown in Fig.~\ref{fig:tpf}, we consider the dynamics of two fluids in the 
domain $\Omega\subset\bR^d$ with 
$\Omega = \Omega_+(t)\cup\Omega_-(t)\cup\Gamma(t)$, where  $d\in\{2,3\}$, 
$\Omega_\pm(t)$ are the regions occupied by the two fluids, and 
$\Gamma(t) = \partial \Omega_-(t)$ is the fluid interface between the two 
fluids. 
Let $\vec u:\Omega\times[0,T]\to\bR^d$ 
be the fluid velocity and $p:\Omega\times[0, T]\to\bR$ be the pressure. The dynamic system is then governed by the standard incompressible Navier--Stokes equations %
\begin{subequations}
\label{eqn:NS}
\begin{alignat}{3}
\label{eq:NS1}
\rho_\pm\partial_t^\bullet\vec u
& = \nabla\cdot\mat{\sigma} + \rho_\pm\,\vec g\qquad &&\mbox{in}\quad\Omega_\pm(t),\\ 
 \nabla\cdot\vec u & = 0\qquad &&\mbox{in}\quad\Omega_\pm(t),
 \label{eq:NS2}
 \end{alignat}
\end{subequations}
where $\rho_\pm$ are the densities of the fluids in $\Omega_\pm(t)$, $\vec g$ is the body acceleration, and $\partial_t^\bullet$ is the material derivative such that for a vector field $\vec\varphi:\Omega\times[0,T]\to\bR^d$
\begin{equation}
\label{eq:MaterD}
\partial_t^\bullet\vec\varphi = \partial_t\vec\varphi + (\vec u\cdot\nabla)\vec\varphi.
\end{equation}
Besides, $\mat{\sigma}$ is the stress tensor  
\begin{equation}
\mat{\sigma}=2\mu_{\pm}\mat{\tD}(\vec u) - p\,\mat{\Id}\quad\mbox{with}\quad\mat{\tD}(\vec u) = \tfrac12\left[\nabla\vec u + (\nabla\vec u)^T\right]\quad\mbox{in}\quad\Omega_\pm(t),
\end{equation}
where $\mu_{\pm}$ are the viscosities of the fluids in $\Omega_\pm(t)$, $\mat{\tD}(\vec u)$ is the strain rate, and $\mat{\Id}\in\bR^{d\times d}$ is the identity matrix.

The fluid interface $\Gamma(t)$ is a hypersurface without boundary, and a parameterization of $\Gamma(t)$ over the reference surface $\Upsilon$  is given by
\begin{equation}
\vecz(\cdot, t): \Upsilon\times[0,T]\to\bR^d.
\end{equation}
Then the induced velocity of the interface is defined as
\begin{equation}
\mathcal{\vec V}(\vecz(\vec q, t ), t) = \partial_t\vecz(\vec q,t)\qquad\mbox{for all}\quad\vec q\in\Upsilon.
\end{equation}
On the fluid interface $\Gamma(t)$, we have 
\begin{align}\label{eq:IFC}
[\vec u]_-^+=\vec 0,\qquad
\bigl[\,\mat{\sigma}\,\vec\nu\bigr]_-^+
= -\gamma\varkappa\,\vec\nu,\qquad
\mathcal{\vv V}\cdot\vec\nu= \vec u\cdot\vec\nu,
\end{align}
where $[\cdot]_-^+$ denotes the jump value from $\Omega_-(t)$ to $\Omega_+(t)$, $\gamma$ and $\varkappa$ are the surface tension and mean curvature of the fluid interface, respectively, and $\vec\nu$ is the unit normal pointing into the region $\Omega_+(t)$. 

Let $\partial\Omega=\partial_1\Omega\cup\partial_2\Omega$ denote the boundary of $\Omega$ with $\partial_1\Omega\cap\partial_2\Omega=\emptyset$. We prescribe a no-slip boundary condition on $\partial_1\Omega$ and a free slip condition on $\partial_2\Omega$ as follows
\begin{subequations}\label{eqn:BD}
\begin{alignat}{3}
\label{eq:BD1}
\vec u &= \vec 0\quad&\mbox{on}\;\partial_{_1}\Omega,\\
\vec u\cdot\vec n=0 ,\quad(\mat{\sigma}\,\vec n)\cdot\vec t&=0\quad\forall\vec t\in\{\vec n\}^\perp\quad&\mbox{on}\;\partial_{_2}\Omega,
\label{eq:BD2}
\end{alignat}
\end{subequations}
where $\vec n$ is the outer unit normal to $\partial\Omega$, and $\{\vec n\}^\perp:=\left\{\vec t\in\bR^d\;:\;\vec t\cdot\vec n=0\right\}$.

Given the initial velocity $\vec u_0 = \vec u(\cdot, 0)$ and the interface $\Gamma_0=\Gamma(0)$, then \eqref{eqn:NS}, together with the interface conditions in \eqref{eq:IFC} and the boundary conditions in \eqref{eqn:BD} form a complete model for the two-phase Navier--Stokes flow. The total free energy of the system consists of the fluid kinetic energy and interface energy 
\begin{equation}
\mathcal{E}(t) = \tfrac{1}{2}\int_\Omega\rho\,|\vec u|^2\,\dL^d + \gamma\int_{\Gamma(t)}\,\dH^{d-1} = \tfrac{1}{2}\int_\Omega\rho\,|\vec u|^2\,\dL^d + \gamma\,|\Gamma(t)|, 
\end{equation}
where we define $\rho(\cdot, t)=\rho_+\mX_{_{\Omega_+(t)}}+\rho_-\mX_{_{\Omega_-(t)}}$ with $\mX_{_E}$ being the usual characteristic function of a set $E$, $\mL^d$ represents the Lebesgue measure in $\bR^d$, and $\mH^{d-1}$ is the $(d-1)$-dimensional Hausdorff measure in $\bR^d$. Moreover, the dynamic system obeys the volume conservation law as well as the energy law
\begin{subequations} 
\begin{align}
\label{eq:VolumeCon}
&\ddt{\vol(\Omega_-(t))}=\int_{\Omega_-(t)}\nabla\cdot\vec u\,\dL^d=0,\\
&\ddt \mathcal{E}(t) = -2\int_\Omega\mu\,\mat{\tD}(\vec u):\mat{\tD}(\vec u)\,\dL^d + \int_\Omega\rho\,\vec u\cdot\vec g\,\dL^d,
\label{eq:EnergyLaw}
\end{align}
\end{subequations}
where $\mu(\cdot, t)=\mu_+\mX_{_{\Omega_+(t)}}+\mu_-\mX_{_{\Omega_-(t)}}$,
see, e.g., \cite{BGN15stable,Barrett20}. 

The main aim of this work is to devise structure-preserving discretizations for
the incompressible two-phase flow problem so that the two physical laws in 
\eqref{eq:VolumeCon} and \eqref{eq:EnergyLaw} are satisfied as well on the 
discrete level. In the following, we consider unfitted and fitted finite 
element approximations in Section~\ref{sec:unfem} and Section~\ref{sec:fitfem},
respectively.

\setcounter{equation}{0}
\section{The unfitted mesh approach}\label{sec:unfem}

\subsection{An Eulerian weak formulation}\label{sec:Eulerweak}

In order to introduce the weak formulation, we define the  following function spaces
\begin{subequations}\label{eqn:UPspaces}
\begin{align}
\mathbb{U}&:=\bigl\{\vec\varphi\in [H^1(\Omega)]^d:\;\vec\varphi = \vec 0\;\;\mbox{on}\;\;\partial_1\Omega,\quad\vec\varphi\cdot\vec n = 0\;\;\mbox{on}\;\;\partial_2\Omega\bigr\},\\
\mathbb{P}&:= \bigl\{\eta\in L^2(\Omega):\;(\eta,~1) = 0\bigr\},\qquad 
\mathbb{V}:=H^1(0,T;[L^2(\Omega)]^d)\cap L^2(0,T;\mathbb{U}),
\end{align}
\end{subequations}
where we denote by $(\cdot,\cdot)$ the $L^2$-inner product over $\Omega$. For all $\vec\chi\in\mathbb{V}$, it holds that (see \cite[(2.16)]{BGN15insoluble})
\begin{equation}
\bigl(\rho\,\partial_t^\bullet\vec u,~\vec\chi\bigr)=\tfrac{1}{2}\left[\ddt\bigl(\rho\,\vec u,~\vec\chi\bigr)+ \bigl(\rho\,\partial_t\vec u,~\vec\chi\bigr) - \bigl(\rho\,\vec u,~\partial_t\vec\chi\bigr)\right] + \mathscr{A}\bigl(\rho,\vec u;~\vec u,\vec\chi\bigr),\label{eq:wim1}
\end{equation}
where we introduced the antisymmetric term 
\begin{equation}
\mathscr{A}(\rho,\vec v;\vec u,\vec\chi) = \tfrac{1}{2}\left[\bigl(\rho\,(\vec v\cdot\nabla)\vec u,~\vec \chi\bigr)-\bigl(\rho\,(\vec v\cdot\nabla)\vec \chi,~\vec u\bigr)\right].\label{eq:antisym}
\end{equation}

We denote by  $\langle\cdot,\cdot\rangle_{\Gamma(t)}$ the $L^2$-inner product over $\Gamma(t)$. Then for the viscous term  in \eqref{eq:NS1}, we take the inner product  with $\vec\chi\in\mathbb{V}$, integrate by parts and  obtain  
\begin{align}
\bigl(\nabla\cdot\mat{\sigma},~\vec\chi\bigr)&=\int_{\Omega_-(t)}(\nabla\cdot\mat{\sigma})\cdot\vec\chi\,\dL^d + \int_{\Omega_+(t)}(\nabla\cdot\mat{\sigma})\cdot\vec\chi\,\dL^d\nn\\
&=-2\bigl(\mu\,\mat{\tD}(\vec u),~\mat{\tD}(\vec\chi)\bigr) + \bigl(p,~\nabla\cdot\vec\chi\bigr) + \gamma\big\langle\varkappa\,\vec\nu,~\vec\chi\big\rangle_{\Gamma(t)},\label{eq:wim2}
\end{align}
on recalling \eqref{eq:IFC} and \eqref{eqn:BD}.

We then propose an {\bf Eulerian weak formulation} for the system of two-phase Navier--Stokes flow as follows. Given $\Gamma(0)=\Gamma_0$ and $\vec u(\cdot,0) = \vec u_0$, we find $\Gamma(t)=\vecz(\Upsilon,t)$ with $\mathcal{\vec V}(\cdot,t)\in[H^1(\Gamma(t)]^d$, $\vec u\in\mathbb{V}$, $p \in L^2(0,T; \mathbb{P})$ and $\varkappa(\cdot,t)\in L^2(\Gamma(t)$ such that (see  \cite[(3.9), (3.10)]{BGN15stable}) for all $t\in(0,T]$
 \begin{subequations}
\label{eqn:weak}
\begin{alignat}{2}
&\tfrac{1}{2}\left[\ddt\bigl(\rho\,\vec u,~\vec\chi\bigr)+ \bigl(\rho\,\partial_t\vec u,~\vec\chi\bigr) - \bigl(\rho\,\vec u,~\partial_t\vec\chi\bigr)\right]+\mathscr{A}\bigl(\rho,\vec u;\vec u,\vec\chi\bigr)\nn\\
&\hspace{1cm}+2\bigl(\mu\,\mat{\tD}(\vec u),~\mat{\tD}(\vec\chi)\bigr)-\bigl(p,~\nabla\cdot\vec\chi\bigr) - \gamma\big\langle\varkappa\,\vec\nu,~\vec\chi\big\rangle_{\Gamma(t)} = \bigl(\rho\,\vec g,~\vec\chi\bigr)\qquad\forall\vec\chi\in\mathbb{V},\label{eq:weak1}\\[0.5em]
&\hspace{3cm}\bigl(\nabla\cdot\vec u,~q\bigr)=0\qquad\forall q\in\mathbb{P},\label{eq:weak2}\\[0.5em]
&\hspace{2.2cm}\big\langle[\mathcal{\vv V}-\vec u]\cdot\vec\nu,~\varphi\big\rangle_{\Gamma(t)}=0\qquad\forall\varphi\in L^2(\Gamma(t)),\label{eq:weak3}\\[0.5em]
&\hspace{0.5cm}\big\langle\varkappa\,\vec\nu,~\vec\zeta\big\rangle_{\Gamma(t))} + \big\langle\nabla_s\vec\id,~\nabla_s\vec\zeta\big\rangle_{\Gamma(t)}=0\qquad\forall\vec\zeta\in [H^1(\Gamma(t))]^d.\label{eq:weak4}
\end{alignat}
\end{subequations}
Here \eqref{eq:weak1} is a direct result from \eqref{eq:wim1} and \eqref{eq:wim2}, and \eqref{eq:weak2} is the divergence free condition \eqref{eq:NS2}. Equation \eqref{eq:weak3} results from the kinetic equation in \eqref{eq:IFC} and \eqref{eq:weak4} is due to the curvature formulation in \eqref{eq:CF}.

It follows from the Reynolds transport theorem that (see \cite{Barrett20})
\begin{subequations} \label{eqn:RT}
\begin{align}\label{eq:tOmega}
&\ddt\vol(\Omega_-(t)) = \big\langle\mathcal{\vv V}\cdot\vec\nu,~1\big\rangle_{\Gamma(t)},\\
&\ddt|\Gamma(t)| = \bigl\langle\nabla_s\vec\id, \nabla_s\mathcal{\vv V}\big\rangle_{\Gamma(t)}.\label{eq:tGamma}
\end{align}
\end{subequations}
Let 
\begin{equation}
\omega(t)=\frac{\vol(\Omega_-(t))}{\vol(\Omega)},\qquad t\in[0,T].\label{eq:lambda}
 \end{equation}
Then it is easy to show that $(\mX_{_{\Omega_-(t)}} - \omega(t))\in\mathbb{P}$. 
Choosing $q = \mX_{_{\Omega_-(t)}} - \omega(t)$ in \eqref{eq:weak2}, $\varphi=1$ in \eqref{eq:weak3} and recalling \eqref{eq:tOmega} yields that \eqref{eq:VolumeCon} is satisfied. Moreover, setting $\vec\chi=\vec u$ in \eqref{eq:weak1}, $q = p$ in \eqref{eq:weak2}, $\varphi = \gamma\varkappa$ in \eqref{eq:weak3} and $\vec\zeta=\mathcal{\vv V}$ in \eqref{eq:weak4}, it is not difficult to show that \eqref{eq:EnergyLaw} is satisfied on noting the identity \eqref{eq:tGamma}.
This implies that the volume conservation and energy law  are satisfied well within the weak formulation.

\subsection{The discretization}\label{sec:ufdis}

We partition the time domain uniformly as $[0,T]=\bigcup_{m=1}^M [t_{m-1},~t_m]$ with $t_m = m\ttau$ and the time step size $\ttau = \frac{T}{M}$. We then solve \eqref{eqn:weak} for $\mathcal{\vv V}$ and $\varkappa$ on the interface $\Gamma$, and for $\vec u$ and $p$ in $\Omega$, using the finite element method.

\vspace{0.2cm}
  {\bf The interface discretization}: Let $\Gamma^m: = \bigcup_{j=1}^{J_\Gamma}\overline{\sigma_j^m}$ be a $(d-1)$-dimensional polyhedral surface to approximate the hypersurface $\Gamma(t_m)$ with vertices $Q_\Gamma^m=\{\vec q_k^m\}_{k=1}^{K_\Gamma}$, where $\{\sigma_j^m\}_{j=1}^{J_\Gamma}$ are mutually disjoint $(d-1)$-simplices in $\bR^d$. On the polyhedral surface $\Gamma^m$, we define the finite element space
\begin{equation}
V^h(\Gamma^m) :=\bigl\{\varphi\in C(\Gamma^m): \varphi|_{\sigma_j^m}\;\;\mbox{is affine}\quad\forall  1\leq j\leq J_\Gamma \bigr\}.
\end{equation}
Then we define the new polyhedral surface 
$\Gamma^{m+1}:=\vec X^{m+1}(\Gamma^m)$ for a $\vec X^{m+1} \in
[V^h(\Gamma^m)]^d$ that is to be determined.
In addition, let
$\left\{\vec q_{j_k}^{m}\right\}_{k=0}^{d-1}$ be the vertices of $\sigma_j^{m}$, ordered with the same orientation for all $\sigma_j^{m}$, $j=1,\ldots, J_\Gamma$. For simplicity, we denote $\sigma_j^{m}=\Delta\left\{\vec q_{j_k}^{m}\right\}_{k=0}^{d-1}$. Then we introduce the unit normal $\vec{\nu}^m$ to $\Gamma^m$; that is,
\begin{equation}\label{eq:vG}
\vec{\nu}^m_{j} := \vec{\nu}^m \mid_{\sigma^{m}_j} :=
\frac{\vec A\{\sigma_j^{m}\}}{
|\vec A\{\sigma_j^{m}\}|}\quad\mbox{ with}\quad \vec A\{\sigma_j^{m}\}=( \vec{q}^{m}_{j_1} - \vec{q}^{m}_{j_0} ) \wedge \ldots \wedge
( \vec{q}^{m}_{j_{d-1}} - \vec{q}^{m}_{j_0}),
\end{equation}
where $\wedge$ is the wedge product and $\vec A\{\sigma_j^{m}\}$ is the orientation vector of $\sigma_j^{m}$.  To approximate the inner product $\langle\cdot,\cdot\rangle_{\Gamma(t_m)}$, we introduce the inner products 
$\langle\cdot,\cdot\rangle_{\Gamma^m}$ and  $\langle\cdot,\cdot\rangle_{\Gamma^m}^h$ over
the current polyhedral surface $\Gamma^m$ via 
\begin{subequations}
\begin{align}
\label{eq:erule}
\langle u, v\rangle_{\Gamma^m} &:=  
\int_{\Gamma^m} u\cdot v \dH^{d-1},\\
\langle u, v \rangle^h_{\Gamma^m} &:=  
\frac{1}{d}\sum_{j=1}^{J_\Gamma} |\sigma^{m}_j|
\sum_{k=0}^{d-1} 
\underset{\sigma^{m}_j\ni \vec{p}\to \vec{q}^{m}_{j_k}}{\lim}\, 
(u\cdot v)(
\vec{p}),
\label{eq:tprule}
\end{align}
\end{subequations}
where $u,v$ are piecewise continuous, with possible jumps
across the edges of $\{\sigma^{m}_j\}_{j=1}^{J_\Gamma}$, and \linebreak
$|\sigma^{m}_j| = \frac{1}{(d-1)!}\,|\vec A\{\sigma_j^{m}\}|$ 
is the measure of $\sigma^{m}_j$.

Following the work in \cite{BZ21SPFEM,BGNZ23}, we introduce a family of polyhedral surfaces $\Gamma^h(t)$ via the linear interpolation between $\Gamma^m$ and $\Gamma^{m+1}$: 
 \begin{equation}
 \Gamma^h(t):=\frac{t_{m+1}-t}{\ttau}\Gamma^m + \frac{t-t_m}{\ttau}\Gamma^{m+1}\quad\mbox{and}\quad \Gamma^h(t)=\cup_{j=1}^{J_\Gamma}\overline{\sigma_j^h}(t),\quad t\in[t_m,~t_{m+1}],\nn
 \end{equation}
where $\{\sigma_j^h(t)\}_{j=1}^{J_\Gamma}$ are mutually disjoint $(d-1)$-simplices and the vertices $\{\vec q_k^{h}(t)\}_{k=1}^{K_\Gamma}$ of $\Gamma^h(t)$ are given by  
\begin{equation}
\label{eq:qh}
\vec q_k^{h}(t) = \frac{t_{m+1} - t}{\Delta t}\vec q_k^{m} + \frac{t - t_m}{\Delta t}\,\vec q_k^{m+1},\quad t\in[t_m,~t_{m+1}],\quad k = 1,\ldots, K_\Gamma.
\end{equation}
 We then define the time-weighted normals  $\vec\nu^{m+\frac{1}{2}}\in  [L^\infty(\Gamma^{m})]^d$ 
\begin{align}
\label{eq:weightv}
\vec\nu^{m+\frac{1}{2}}|_{\sigma_j^{m}}=\vec\nu_{j}^{m+\frac{1}{2}}&:=\frac{1}{\Delta t\,|\vec A\{\sigma_j^{m}\}|}\int_{t_{m}}^{t_{m+1}}\vec A\{\sigma_j^{h}(t)\}\,\rd t\qquad \forall 1\leq j \leq J_\Gamma.
\end{align}
Let $\Omega_-^m$ be the interior of $\Gamma^m$ and $\Omega_+^m$ be the exterior of $\Gamma^m$ in $\Omega$. Then we have the following lemma.
\begin{lem}\label{lem:vc}
Let $\vec X^{m+1}\in [V^h(\Gamma^{m})]^d$. 
Then it holds
\begin{subequations}
\begin{equation}
\label{eq:Vlem}
\big\langle(\vec X^{m+1} - \vec\id)\cdot\vec\nu^{m+\frac{1}{2}},~1\big\rangle_{\Gamma^{m}}^h=\vol(\Omega_-^{m+1})-\vol(\Omega_-^{m}).
\end{equation}
Moreover, it holds that
\begin{equation}
\big\langle\nabla_s\vec X^{m+1},~\nabla_s(\vec X^{m+1}-\vec\id)\big\rangle_{\Gamma^m}\geq |\Gamma^{m+1}| - |\Gamma^m|.\label{eq:Elem}
\end{equation}
\end{subequations}
\end{lem}
Here \eqref{eq:Vlem} can be regarded as a natural discrete analogue of \eqref{eq:tOmega}, and its proof can be found in \cite{BZ21SPFEM,mullins}. While for the stability bound in \eqref{eq:Elem}, the reader can refer to \cite[Lemma 1]{Bansch01} and \cite[Lemma 57]{Barrett20} for the proof.  

\vspace{0.2cm}
 {\bf The bulk discretization}: We next consider the partition of the bulk domain. At time $t_m$, a regular partition of $\Omega$ with $K^m_\Omega$ vertices is given by
\begin{equation}
\overline{\Omega} = \cup_{j = 1}^{J_\Omega^m} \overline{o_j^m}\quad\mbox{with}\quad Q^m=\{\vec a_k^m\}_{k=1}^{K^m_\Omega},
\end{equation}
where $\{o_j^m\}_{j=1}^{J_\Omega^m}$  are the mutually disjoint open simplices in $\bR^d$. We define the bulk mesh $\mathscr{T}^m=\bigl\{o_j^m: j = 1,\ldots,J_\Omega^m\bigr\}$, and then introduce the finite element spaces associated with  $\mathscr{T}^m$ as
\begin{subequations}
\begin{align}
S_k^m &:= \bigl\{\chi\in C(\overline{\Omega}):\quad \chi|_{o_j^m}\in P_k(o_j^m)\quad\forall 1\leq j\leq J_\Omega^m\bigr\}, \quad k\in\bN_+,\\
S_0^m &:= \bigl\{\chi\in L^2(\Omega):\quad \chi|_{o_j^m}\;\;\mbox{is constant}\quad\forall 1\leq j \leq J_\Omega^m\bigr\},
\end{align}
\end{subequations}
where $P_k(o_j^m)$ denotes the space of polynomials of degree at most $k$ on $o_j^m$. Let $\mathbb{U}^m$ and $\mathbb{P}^m$ be the finite element spaces for the numerical solutions of the fluid velocity and pressure,
recall \eqref{eqn:UPspaces}. 
It is natural to consider the following pair elements \cite{Agnese16,Agnese20} 
\begin{subequations}\label{eq:UPspaces}
\begin{align}\label{eq:P2P1}
\mbox{P2-P1}:&\quad\bigl(\mathbb{U}^m,~\mathbb{P}^m\bigr)= \bigl([S_2^m]^d\cap\mathbb{U},~S_1^m\cap\mathbb{P}\bigr),\\\label{eq:P2P0}
\mbox{P2-P0}:&\quad\bigl(\mathbb{U}^m,~\mathbb{P}^m\bigr)= \bigl([S_2^m]^d\cap\mathbb{U},~S_0^m\cap\mathbb{P}\bigr),\\\label{eq:P2P1P0}
\mbox{P2-(P1+P0)}:&\quad\bigl(\mathbb{U}^m,~\mathbb{P}^m\bigr)= \bigl([S_2^m]^d\cap\mathbb{U},~(S_1^m+S_0^m)\cap\mathbb{P}\bigr),
\end{align}
which usually satisfy the LBB inf-sup stability condition. 
In particular, for the LBB stability of (\ref{eq:P2P1}) 
we refer to \cite[p.\ 252]{BrezziF91} for 
$d=2$ and to \cite{Boffi97} for $d=3$, 
while the stability of (\ref{eq:P2P0}) is shown in 
\cite[p.\ 221]{BrezziF91} for $d=2$. 
The LBB stability of (\ref{eq:P2P1P0}) is shown in 
\cite{Boffi2012local} for $d=2$.
Here the results for (\ref{eq:P2P1},c) need the weak
constraint that all the elements of the bulk mesh have a vertex in $\Omega$. Later in the unfitted approximations  we actually need to enrich the pressure space with additional degrees of freedom, see \eqref{eq:Asup2} and  Remark \ref{rem:ufXFEM}. In particular, a possible choice for the enriched elements is 
\begin{align}
\mbox{P2-P1 with XFEM}:&\quad\bigl(\mathbb{U}^m,~\mathbb{P}^m\bigr)= \bigl([S_2^m]^d\cap\mathbb{U},~\mbox{span}(S_1^m\cup\{\mX_{_{\Omega_-^m}}\})\cap\mathbb{P}\bigr).\label{eq:P2P1XFEM}
\end{align}
\end{subequations}

\subsection{The Eulerian structure-preserving method }
In the unfitted approach, the bulk mesh $\mathscr{T}^m$ is decoupled from the polyhedral surface $\Gamma^m$. We partition the elements of the bulk mesh $\mathscr{T}^m$ into interior, exterior and interface elements as
 \begin{equation}
 \mathscr{T}_-^m:=\{o\in\mathscr{T}^m: o\subset\Omega_-^m\},\qquad \mathscr{T}_+^m=\{o\in\mathscr{T}^m: o\subset\Omega_+^m\},\qquad \mathscr{T}_\Gamma^m=\{o\in\mathscr{T}^m: o\cap\Gamma^m\neq\emptyset\},\nn
 \end{equation}
where $\Omega_-^m$ and $\Omega_+^m$ denote the interior and the exterior of
$\Gamma^m$, respectively.
Let $\rho^m$ and $\mu^m$ be the numerical approximations of the density 
$\rho(\cdot, t)$ and viscosity $\mu(\cdot, t)$ at $t=t_m$, respectively. 
Then we define $\rho^m\in S_0^m$, $\mu^m\in S_0^m$ such that
\begin{equation}
\label{eqn:discreterhovis}
\rho^m|_{o}:=\left\{
\begin{array}{ll}
\rho_-, & \text{if}\ o\in \mathscr{T}_-^m, \vspace{0.15cm}\\
\rho_+, & \text{if}\ o\in\mathscr{T}_+^m,\vspace{0.15cm}\\
\frac{1}{2}(\rho_-+\rho_+), &\text{if}\ o\in \mathscr{T}_{\Gamma}^m,
\end{array}\right.\quad\mbox{and}\quad
\mu^m|_{o}:=\left\{\begin{array}{ll}
\mu_-, &\text{if}\ o\in \mathscr{T}_1^m,\vspace{0.15cm}\\
\mu_+,&\text{if}\ o\in\mathscr{T}_2^m, \vspace{0.15cm}\\
\frac{1}{2}(\mu_-+\mu_+), & \text{if}\ o\in \mathscr{T}_{\Gamma}^m.
\end{array}\right.\nonumber
\end{equation} 

We introduce the standard interpolation operator $I_k^m: C(\overline{\Omega})\to [S_k^m]^d$ for $k\geq 1$, and the standard projection operator $I_0^m: L^1(\Omega)\to S_0^m$ with $(I_0^m\eta)|_o= \frac{1}{|o|}\int_o\eta\dL^d$ for $o\in\mathscr{T}^m$. Then the unfitted finite element approximation of \eqref{eqn:weak} is given as follows. 
Let $\vec U^0\in\mathbb{U}^0$ be an approximation of the initial fluid 
velocity $\vec u_0$. 
Moreover, let $\Gamma^0$ be a polyhedral approximation of the initial fluid 
interface $\Gamma_0$ and set $\vec X^0 = \vec\id_{| \Gamma^0}$. We also 
set $\rho^{-1} = \rho^0$. Then, for $m \geq 0$, find $U^{m+1}\in\mathbb{U}^m$, $P^{m+1}\in\mathbb{P}^m$, $\vec X^{m+1}\in [V^h(\Gamma^m)]^d$ and $\kappa^{m+1}\in V^h(\Gamma^m)$ such that
\begin{subequations}\label{eqn:fem}
\begin{align}
&\tfrac{1}{2}\Bigl[\Bigl(\frac{\rho^m\vec U^{m+1}-I_0^m\rho^{m-1}
I_2^m\vec U^m}{\ttau}+I_0^m\rho^{m-1}\frac{\vec U^{m+1}-I_2^m\vec U^m}
{\ttau},~\vec\chi^h\Bigr)\Bigr]+\mathscr{A}(\rho^m, I_2^m\vec U^{m};~\vec U^{m+1},\vec\chi^h\big)\nn\\ 
&\hspace{1.5cm}+2\bigl(\mu^{m}\mat{\tD}(\vec U^{m+1}),~\mat{\tD}(\vec\chi^h)\bigr) - \bigl(P^{m+1},~\nabla\cdot\vec\chi^h\bigr)- \gamma\big\langle\kappa^{m+1}\,\vec\nu^m,~\vec\chi^h\big\rangle_{\Gamma^m}\nn\\&\hspace{3cm} = \bigl(\rho^m\vec g,~\vec\chi^h\bigr)\qquad\forall\vec\chi^h\in\mathbb{U}^m,\label{eq:fem1}\\[0.5em]
&\hspace{0.5cm}\bigl(\nabla\cdot\vec U^{m+1},~q^h\bigr) = 0\qquad \forall q^h \in \mathbb{P}^m, \label{eq:fem2}\\[0.5em]
&\hspace{0.5cm}\big\langle\frac{\vec X^{m+1}-\vec\id}{\ttau}\cdot\vec\nu^{m+\frac{1}{2}},~\varphi^h\big\rangle_{\Gamma^m}^h - \big\langle\vec U^{m+1}\cdot\vec\nu^m,~\varphi^h\big\rangle_{\Gamma^m} =0\qquad\forall\varphi^h\in V^h(\Gamma^m),\label{eq:fem3}\\[0.5em]
&\hspace{0.5cm}\big\langle\kappa^{m+1}\,\vec\nu^{m+\frac{1}{2}},~\vec\zeta^h\big\rangle_{\Gamma^m}^h + \big\langle\nabla_s\vec X^{m+1},~\nabla_s\vec\zeta^h\big\rangle_{\Gamma^m} =0\qquad\forall\vec\zeta^h\in [V^h(\Gamma^m)]^d,\label{eq:fem4}
\end{align}
\end{subequations}
and then we set $\Gamma^{m+1}:=\vec X^{m+1}(\Gamma^m)$. 
The above numerical method \eqref{eqn:fem} is very similar to the BGN method 
in \cite{BGN15stable}, except that here the time-weighted discrete normals 
$\vec\nu^{m+\frac12}$ are employed in the first terms of \eqref{eq:fem3} and 
\eqref{eq:fem4}, instead of the explicit discrete normals $\vec\nu^m$
in \cite{BGN15stable}. 
Therefore, the scheme \eqref{eqn:fem} not only inherits the unconditional
stability property of the original BGN method, but also has the property of 
exact volume preservation, as shown by the following theorem.

\begin{thm}[stability and volume conservation]\label{thm:stavc}
Let $(\vec U^{m+1},P^{m+1},\vec X^{m+1}, \kappa^{m+1})$ be a solution of \eqref{eqn:fem} for $m=0,1,\ldots, M-1$. Then it holds 
\begin{equation}
\mathcal{E}(\rho^m,\vec U^{m+1}, \Gamma^{m+1}) + 2\ttau \norm{\sqrt{\mu^m}\mat{\tD}(\vec U^{m+1})}\leq \mathcal{E}(I_0^m\rho^{m-1}, I_2^m\vec U^m, \Gamma^m) + \ttau\left(\rho^m\,\vec g,~\vec U^{m+1}\right), 
\end{equation}
where $\mathcal{E}(\rho,\vec U, \Gamma) = \frac{1}{2}\norm{\sqrt{\rho}\,\vec U} + \gamma|\Gamma|$ with $\norm{\cdot}$ being the norm induced by 
the inner product $(\cdot,\cdot)$. In addition, on assuming that 
\begin{equation}
\mathcal{E}(I_0^m\rho^{m-1}, I_2^m\vec U^m, \Gamma^m) \leq \mathcal{E}(\rho^{m-1}, \vec U^m, \Gamma^m)\quad\mbox{for}\quad m = 1,\ldots, M-1,\label{eq:ufassup}
\end{equation}
it holds that
\begin{align}
\mathcal{E}(\rho^k,\vec U^{k+1},\Gamma^{k+1}) + 2\ttau\,\sum_{m=0}^k\norm{\sqrt{\mu^m}\mat{\tD}(\vec U^{m+1})} \leq \mathcal{E}(\rho^{0},\vec U^0,\Gamma^0) + \ttau\,\sum_{m=0}^k\bigl(\rho^{m}\,\vec g, ~\vec U^{m+1}\bigr),\label{eq:gDec}
\end{align}
for $k=0,\ldots,M-1$. Moreover, it holds for $m = 0,\ldots, M-1$ that if
\begin{equation}
\left(\mX_{_{\Omega_-^m}}-\omega^m\right)\in\mathbb{P}^m\qquad{\rm with}\quad \omega^m =\frac{\vol(\Omega_-^m)}{\vol(\Omega)},\label{eq:Asup2}
\end{equation}
then
\begin{equation}
\vol(\Omega_-^{m+1}) = \vol(\Omega_-^m).
\label{eq:Dvc}
\end{equation}
\end{thm}
\begin{proof}
Choosing $\vec\chi^h = \ttau\,\vec U^{m+1}$ in \eqref{eq:fem1}, $q^h = P^{m+1}$ in \eqref{eq:fem2}, $\varphi^h = \ttau\,\gamma\kappa^{m+1}$ in \eqref{eq:fem3} and $\vec\zeta^h = (\vec X^{m+1} - \vec\id_{|\Gamma^m})$ in \eqref{eq:fem4}, and combining these equations yields 
\begin{align}
&\frac{1}{2}\bigl(\rho^m\vec U^{m+1}-I_0^m\rho^{m-1}
I_2^m\vec U^m+I_0^m\rho^{m-1}[\vec U^{m+1}-I_2^m\vec U^m],~\vec U^{m+1}\bigr) + 2\ttau\,\bigl(\mu^{m}\,\mat{\tD}(\vec U^{m+1}),~\mat{\tD}(\vec U^{m+1})\bigr)\nn\\
&\hspace{2cm}+\gamma\big\langle\nabla_s\vec X^{m+1},~\nabla_s(\vec X^{m+1}-\vec\id)\big\rangle_{\Gamma^m} = \ttau\,\bigl(\rho^m\,\vec g,~\vec U^{m+1}\bigr).\label{eq:E1}
\end{align}
It is easy to show that
\begin{align}
&\bigl(\rho^m\vec U^{m+1}-I_0^m\rho^{m-1}
I_2^m\vec U^m+I_0^m\rho^{m-1}(\vec U^{m+1}-I_2^m\vec U^m),~\vec U^{m+1}\bigr)\nn \\&= \bigl(\rho^m\,\vec U^{m+1},~\vec U^{m+1}\bigr) - \bigl(I_0^m\rho^{m-1}\,I_2^m\vec U^m,~I_2^m\vec U^m\bigr) + \bigl(I_0^m\rho^{m-1}[\vec U^{m+1}-I_2^m\vec U^m],~[\vec U^{m+1}-I_2^m\vec U^m]\bigr)\nn\\
&\geq \bigl(\rho^m\,\vec U^{m+1},~\vec U^{m+1}\bigr) - \bigl(I_0^m\rho^{m-1}\,I_2^m\vec U^m,~I_2^m\vec U^m\bigr).\label{eq:E2}
\end{align}
Using \eqref{eq:E2} and \eqref{eq:Elem} in \eqref{eq:E1}, we immediately obtain 
\begin{equation}
\mathcal{E}(\rho^m,\vec U^{m+1}, \Gamma^{m+1}) + 2\ttau\norm{\sqrt{\mu^m}\,\mat{\tD}(\vec U^{m+1})} \leq \mathcal{E}(I_0^m\rho^{m-1}, I_2^m\vec U^m, \Gamma^m) + \ttau\bigl(\rho^m\,\vec g,~\vec U^{m+1}\bigr),
\label{eq:Des}
\end{equation}
Summing up for $m=0,\ldots, k$ gives \eqref{eq:gDec} immediately on recalling 
\eqref{eq:ufassup} and noting that \eqref{eq:ufassup} trivially holds for
$m=0$. 

We next choose $q^h = \mX_{_{\Omega_-^m}}-\omega^m$ in \eqref{eq:fem2} and obtain
\begin{equation}
0=\bigl(\nabla\cdot\vec U^{m+1},~\mX_{_{\Omega_-^m}}-\omega^m\bigr) = \int_{\Omega_-^m}\nabla\cdot\vec U^{m+1}\,\dL^d - \omega^m\int_\Omega\nabla\cdot\vec U^{m+1}\,\dL^d,\nn
\end{equation}
which then yields
\begin{equation}
0=\big\langle\vec U^{m+1}\cdot\vec\nu^m,~1\big\rangle_{\Gamma^m}.\label{eq:V1}
\end{equation}
On the other hand setting $\varphi^h = \ttau$ in \eqref{eq:fem3} yields
\begin{align}
\big\langle(\vec X^{m+1}-\vec\id)\cdot\vec\nu^{m+\frac{1}{2}},~1\big\rangle_{\Gamma^m}^h = \ttau\big\langle\vec U^{m+1}\cdot\vec\nu^m,~1 \big\rangle_{\Gamma^m},\nn
\end{align} 
which implies  \eqref{eq:Dvc} on recalling \eqref{eq:V1} and \eqref{eq:Vlem}.
\end{proof}

\begin{rem} It is desirable to achieve the global stability bound in \eqref{eq:gDec}, which requires the assumption \eqref{eq:ufassup}. In fact, \eqref{eq:ufassup} is well satisfied in the following cases
\begin{itemize}
\item[i)] No mesh adaptation is employed, i.e., $\mathscr{T}^m=\mathscr{T}^0$ for $m=1,\ldots,M-1$.
\item[ii)] Mesh refinement routines without coarsening are employed so that $\mathbb{U}^{m-1}\subset\mathbb{U}^m$.
\end{itemize}
Nevertheless, in practical computation, we implement a bulk mesh adaptation strategy, which is described in detail in \cite{BGN15stable}. 
\end{rem}

\begin{rem}\label{rem:ufXFEM}
The exact volume preservation requires
\begin{itemize}
\item [(i)] the exact evaluation of the time-weighted normals in \eqref{eq:weightv};
\item [(ii)] the assumption \eqref{eq:Asup2}.
\end{itemize}
Condition (i) is guaranteed by noting that $\vec A\{\sigma_j^h(t)\}$ is a polynomial of degree $d-1$ for the variable $t$, recall \eqref{eq:vG}, and thus the integration can be calculated exactly via appropriate quadrature rules. To satisfy condition (ii), we employ an XFEM procedure, using the finite element spaces in \eqref{eq:P2P1XFEM}, where the pressure approximating space is enriched with a single basis function. The contribution of this new basis function to \eqref{eq:fem1} and \eqref{eq:fem2} can be written in terms of integrals over $\Gamma^m$  (see \cite{BGN15stable})
\begin{equation}
\bigl(\nabla\cdot\vec\chi^h,~\mX_{_{\Omega_-^m}}\bigr) = \int_{\Omega_-^m}\nabla\cdot\vec\chi^h\,\dL^d = \int_{\Gamma^m}\vec\chi^h\cdot\vec\nu^m\,\dH^{d-1}\quad\mbox{for all}\quad\vec\chi^h\in\mathbb{U}^m.\nn
\end{equation}
\end{rem}

\begin{rem}
The introduced method \eqref{eqn:fem} is weakly nonlinear and gives rise to a system of nonlinear polynomial equations, which can be solved efficiently via  Picard-type iterative method as follows. Upon setting 
$\vec X^{m+1,0} = \vec\id_{| \Gamma^m}$,  
for each $l\geq 0$, we let $\Gamma^{m+1,l} = \vec X^{m+1,l}(\Gamma^m)$, 
and define $\vec\nu^{m+\frac{1}{2},l}$ through the formula \eqref{eq:weightv} 
with $\Gamma^{m+1}$ replaced by $\Gamma^{m+1,l}$. 
Then we find $\bigl(\vec U^{m+1,l+1}, P^{m+1,l+1}, \vec X^{m+1,l+1}, \kappa^{m+1,l+1}\bigr)\in\mathbb{U}^m\times\mathbb{P}^m\times [V^h(\Gamma^m)]^d\times V^h(\Gamma^m)$ such that
\begin{subequations}\label{eqn:Picardfem}
\begin{align}
&\tfrac{1}{2}\Bigl[\Bigl(\frac{\rho^m\vec U^{m+1, l+1}-I_0^m\rho^{m-1}
I_2^m\vec U^m}{\ttau}+I_0^m\rho^{m-1}\frac{\vec U^{m+1, l+1}-I_2^m\vec U^m}
{\ttau},~\vec\chi^h\Bigr)\Bigr]+\mathscr{A}(\rho^m, I_2^m\vec U^{m};~\vec U^{m+1, l+1},\vec\chi^h\big)\nn\\
&\hspace{1.5cm}+2\bigl(\mu^{m}\,\mat{\tD}(\vec U^{m+1, l+1}),~\mat{\tD}(\vec\chi^h)\bigr) - \bigl(P^{m+1, l+1},~\nabla\cdot\vec\chi^h\bigr)- \gamma\big\langle\kappa^{m+1, l+1}\,\vec\nu^m,~\vec\chi^h\big\rangle_{\Gamma^m}\nn\\&\hspace{2.5cm} = \bigl(\rho^m\vec g,~\vec\chi^h\bigr)\qquad\forall\vec\chi^h\in\mathbb{U}^m,\label{eq:pfem1}\\[0.5em]
&\hspace{0.5cm}\bigl(\nabla\cdot\vec U^{m+1,l+1},~q^h\bigr) = 0\qquad\forall q^h \in \mathbb{P}^m, \label{eq:pfem2}\\[0.5em]
&\hspace{0.5cm}\big\langle\frac{\vec X^{m+1, l+1}-\vec\id}{\ttau}\cdot\vec\nu^{m+\frac{1}{2}, l},~\varphi^h\big\rangle_{\Gamma^m}^h - \big\langle\vec U^{m+1, l+1}\cdot\vec\nu^m,~\varphi^h\big\rangle_{\Gamma^m} =0\qquad\forall\varphi^h\in V^h(\Gamma^m),\label{eq:pfem3}\\[0.5em]
&\hspace{0.5cm}\big\langle\kappa^{m+1, l+1}\,\vec\nu^{m+\frac{1}{2}, l},~\vec\zeta^h\big\rangle_{\Gamma^m}^h + \big\langle\nabla_s\vec X^{m+1, l+1},~\nabla_s\vec\zeta^h\big\rangle_{\Gamma^m} =0\qquad\forall\vec\zeta^h\in [V^h(\Gamma^m)]^d.\label{eq:pfem4}
\end{align}
\end{subequations}
In practical computations, we perform the 
above iteration until 
\begin{equation}
\max_{\vec q\in Q_\Gamma^m}|\vec X^{m+1,l+1}(\vec q)-\vec X^{m+1,l}(\vec q)|\leq {\rm tol}\quad\mbox{and}\quad \norm{\vec U^{m+1,l+1} - \vec U^{m+1,l}}_{\infty}\leq {\rm tol},
\label{eq:iterativeTol}
\end{equation}
where $\norm{\cdot}_{\infty}$ is the $L^{\infty}$-norm on $\Omega$ and ${\rm tol}$ is a chosen tolerance.
\end{rem}

\setcounter{equation}{0}
\section{The fitted mesh approach}\label{sec:fitfem}

\subsection{An ALE weak formulation}\label{sec:ALEweak}

Let $\mathcal{O}\subset\bR^d$ be the ALE reference domain of $\Omega$ and let $\bigl\{\vec{\mathcal{A}}[t]\,\bigr\}_{t\in[0,T]}$ be the family of ALE mappings 
\begin{equation}
\vec{\mathcal{A}}[t]: \mathcal{O}\to \Omega, \qquad \vec y \mapsto \vec{\mathcal{A}}[t](\vec y) = \vec x(\vec y, t)\quad \mbox{for all}\quad t\in[0,T],\quad \vec y\in\mathcal{O},
\end{equation}
with $\Upsilon\subset\mathcal{O}$ such that $\Gamma(t)=\vecz(\Upsilon,t)$. We assume the introduced mappings satisfy $\vec{\mathcal{A}}[t]\in [W^{1,\infty}(\mathcal{O})]^d$, and $\vec{\mathcal{A}}[t]^{-1}\in [W^{1,\infty}(\Omega)]^d$. The domain mesh velocity is defined by
\begin{equation}
\vec w(\vec x, t) :=\left.\frac{\partial\vec x(\vec y,t)}{\partial t} \right|_{\vec y = \vec{\mathcal{A}}[t]^{-1}(\vec x)}\quad\mbox{for all}\quad t\in[0,T]\quad\mbox{and}\quad \vec x\in\Omega.\label{eq:ALEmeshv}
\end{equation}
We allow the ALE mappings to be
somehow arbitrary, except that on the boundary the domain mesh velocity should satisfy 
\begin{equation}
[\vec w(\vec x, t) - \vec u(\vec x, t)]\cdot\vec\nu = 0\quad\mbox{on}\quad \Gamma(t);\qquad [\vec w(\vec x,t)-\vec u(\vec x, t)]\cdot\vec n = 0\quad\mbox{on}\quad\partial\Omega.\label{eq:MeshVBD}
\end{equation}
The construction of the ALE mappings will be presented in  Section \ref{sec:DALE}. Given $\vec\varphi : \Omega \times [0,T] \to \bR^d$, we introduce the time derivative with respect to the ALE reference domain as
\begin{align}
\partial_t^\circ\vec\varphi :=\partial_t \vec\varphi+ (\vec w\cdot\nabla) \vec\varphi.
\label{eq:ALED}
\end{align}
On recalling \eqref{eq:MaterD}, this then yields that
\begin{equation}
\partial_t^\bullet\vec\varphi = \partial_t^\circ\vec\varphi + ([\vec u-\vec w]\cdot\nabla)\vec\varphi.
\label{eq:AMD}
\end{equation}

Denote by $(\cdot,\cdot)_{\Omega_\pm(t)}$ the $L^2$-inner product on $\Omega_\pm(t)$, respectively. 
 By \eqref{eq:AMD}, we have for the inertia term in \eqref{eq:NS1} that
\begin{equation}
\bigl(\partial_t^\bullet\vec u,~\vec\chi\big)_{\Omega_\pm(t)} = \bigl(\partial_t^\circ\vec u,~\vec\chi\bigr)_{\Omega_\pm(t)} + \bigl(([\vec u-\vec w]\cdot\nabla)\vec u,~\vec\chi\bigr)_{\Omega_\pm(t)}\qquad\forall \vec\chi\in [H^1(\Omega)]^d,
\label{eq:nALE1}
\end{equation}
and the second term of \eqref{eq:nALE1} can be reformulated as
\begin{align}
\bigl(([\vec u-\vec w]\cdot\nabla)\vec u,~\vec\chi\bigr)_{\Omega_\pm(t)} &= \tfrac{1}{2}\left[\bigl(([\vec u-\vec w]\cdot\nabla)\vec u,~\vec \chi\bigr)_{\Omega_\pm(t)}-\bigl(([\vec u-\vec w]\cdot\nabla)\vec\chi,~\vec u\bigr)_{\Omega_\pm(t)}\right]\nn\\&\quad+\tfrac{1}{2}\bigl(\vec u - \vec w,~\nabla(\vec u\cdot\vec\chi)\bigr)_{\Omega_\pm(t)}\nn\\
&=\tfrac{1}{2}\left[\bigl(([\vec u-\vec w]\cdot\nabla)\vec u,~\vec \chi\bigr)_{\Omega_\pm(t)}-\bigl(([\vec u-\vec w]\cdot\nabla)\vec\chi,~\vec u\bigr)_{\Omega_\pm(t)}\right] + \tfrac{1}{2}\bigl(\nabla\cdot\vec w,~\vec u\cdot\vec \chi\bigr)_{\Omega_\pm(t)},\nn
\end{align}
where we applied integration by parts and used the divergence free condition \eqref{eq:NS2} as well as \eqref{eq:MeshVBD}. Therefore, we can recast \eqref{eq:nALE1} as
\begin{align}
\bigl(\partial_t^\bullet\vec u,~\vec\chi\bigr)_{\Omega_\pm(t)}& = \bigl(\partial_t^\circ\vec u,~\vec\chi\bigr)_{\Omega_\pm(t)}+ \tfrac{1}{2}\bigl(\nabla\cdot\vec w,~\vec u\cdot\vec\chi\bigr)_{\Omega_\pm(t)}\nn\\&\quad  + \tfrac{1}{2}\left[\bigl(([\vec u-\vec w]\cdot\nabla)\vec u,~\vec \chi\bigr)_{\Omega_\pm(t)}-\bigl(([\vec u-\vec w]\cdot\nabla)\vec\chi,~\vec u\bigr)_{\Omega_\pm(t)}\right].
\label{eq:nALEn}
\end{align}
Multiplying \eqref{eq:nALEn} with $\rho_\pm$ 
and then summing the two equations yields
\begin{equation}
\bigl(\rho\,\partial_t^\bullet\vec u,~\vec\chi\bigr) = \bigl(\rho\,\partial_t^\circ\vec u,~\vec\chi\bigr) + \mathscr{A}\bigl(\rho,\vec u-\vec w; ~\vec u,\vec\chi\bigr) + \tfrac{1}{2}\bigl(\rho\,\nabla\cdot\vec w,~\vec u\cdot\vec\chi\bigr)
\qquad\forall \vec\chi\in [H^1(\Omega)]^d,
\label{eq:nALE}
\end{equation}
where $\mathscr{A}(\rho,\vec v;\vec u,\vec\chi)$ is the antisymmetric term defined in \eqref{eq:antisym}.

Collecting the results in \eqref{eq:nALE} and \eqref{eq:wim2}, it is natural to consider the following {\bf nonconservative ALE weak formulation} for two-phase flow. Given $\Gamma(0)=\Gamma_0$ and $\vec u(\cdot,0) = \vec u_0$, we find $\Gamma(t)=\vecz(\Upsilon,t)$ with $\mathcal{\vec V}(\cdot,t)\in  [H^1(\Gamma(t)]^d$, $\vec u\in\mathbb{V}$, $p \in L^2(0,T; \mathbb{P})$ and $\varkappa(\cdot,t)\in L^2(\Gamma(t)$ such that for all $t\in(0,T]$
\begin{subequations}\label{eqn:ALEweak}
\begin{alignat}{2}
&\bigl(\rho\,\partial_t^\circ\vec u,~\vec\chi\bigr)+\tfrac{1}{2}\bigl(\rho\,\nabla\cdot\vec w,~\vec u\cdot\vec\chi\bigr)+\mathscr{A}\bigl(\rho,\vec u-\vec w;\vec u,\vec\chi\bigr)+2\bigl(\mu\,\mat{\tD}(\vec u),~\mat{\tD}(\vec\chi)\bigr)\nn\\
&\hspace{1cm}-\bigl(p,~\nabla\cdot\vec\chi\bigr) - \gamma\big\langle\varkappa\,\vec\nu,~\vec\chi\big\rangle_{\Gamma(t)} = \bigl(\rho\,\vec g,~\vec\chi\bigr)\qquad\forall\vec\chi\in\mathbb{V},\label{eq:ALEweak1}\\[0.5em]
&\hspace{3cm}\bigl(\nabla\cdot\vec u,~q\bigr)=0\qquad\forall q\in\mathbb{P},\label{eq:ALEweak2}\\[0.5em]
&\hspace{2.2cm}\big\langle[\mathcal{\vv V}-\vec u]\cdot\vec\nu,~\varphi\big\rangle_{\Gamma(t)}=0\qquad\forall\varphi\in L^2(\Gamma(t)),\label{eq:ALEweak3}\\[0.5em]
&\hspace{0.5cm}\big\langle\varkappa\,\vec\nu,~\vec\zeta\big\rangle_{\Gamma(t))} + \big\langle\nabla_s\vec\id,~\nabla_s\vec\zeta\big\rangle_{\Gamma(t)}=0\qquad\forall\vec\zeta\in [H^1(\Gamma(t))]^d.\label{eq:ALEweak4}
\end{alignat}
\end{subequations}
Choosing $q = \mX_{_{\Omega_-(t)}} - \omega(t)\in\mathbb{P}$ in \eqref{eq:ALEweak2} and $\varphi=1\in L^2(\Gamma(t))$ in \eqref{eq:ALEweak3} yields that \eqref{eq:VolumeCon} is satisfied. 
In addition, it follows from \eqref{eq:ALERey} that
\begin{equation}
\tfrac{1}{2}\ddt\bigl(|\vec u|^2,~1\bigr)_{\Omega_\pm(t)}=\bigl(\partial_t^\circ\vec u,~\vec u\bigr)_{\Omega_\pm(t)} + \tfrac{1}{2}\bigl(|\vec u|^2,~\nabla\cdot\vec w\bigr)_{\Omega_\pm(t)}.
\label{eq:DK1}
\end{equation}
Multiplying \eqref{eq:DK1} with $\rho_\pm$ and summing the two equations yields 
\begin{equation}
\tfrac{1}{2}\ddt\bigl(\rho\,|\vec u|^2,~1\bigr) = \bigl(\rho\,\partial_t^\circ\vec u,~\vec u\bigr) + \tfrac{1}{2}\bigl(\rho\,|\vec u|^2,~\nabla\cdot\vec w\bigr).\label{eq:tKinetic}
\end{equation}
Now choosing $\vec \chi = \vec u$ in \eqref{eq:ALEweak1}, $q = p$ in \eqref{eq:ALEweak2}, $\varphi = \gamma\varkappa$ in \eqref{eq:ALEweak3} and $\vec\zeta=\mathcal{\vv V}$ in \eqref{eq:ALEweak4}, and noting \eqref{eq:tKinetic}
and \eqref{eq:tGamma}, we see that \eqref{eq:EnergyLaw} holds as well.

\begin{rem}\label{rem:DL}
In the case that $\nabla\cdot\vec w = 0$, the second term in \eqref{eq:ALEweak1} disappears. Then the above ALE weak formulation \eqref{eqn:ALEweak} will collapse to the formulation proposed in \cite{Duan2022energy}.
\end{rem}

\begin{rem}\label{rem:cALE}
In view of the work for a single fluid phase in \cite{ivanvcic2022energy}, it is also possible to consider an alternative ALE weak formulation that is in conservative form. Precisely the formulation is given by \eqref{eqn:ALEweak} except that we replace  \eqref{eq:ALEweak1} with
\begin{align}
&\ddt\bigl(\rho\,\vec u,~\vec\chi\bigr)-\tfrac{1}{2}\bigl(\rho\,\nabla\cdot\vec w,~\vec u\cdot\vec\chi\bigr) - \bigl(\rho\,\vec u,~\partial_t^\circ\vec\chi\bigr)+\mathscr{A}\bigl(\rho,\vec u-\vec w;\vec u,\vec\chi\bigr)+2\bigl(\mu\,\mat{\tD}(\vec u),~\mat{\tD}(\vec\chi)\bigr)\nn\\
&\hspace{1cm}-\bigl(p,~\nabla\cdot\vec\chi\bigr) - \gamma\big\langle\varkappa\,\vec\nu,~\vec\chi\big\rangle_{\Gamma(t)} = \bigl(\rho\,\vec g,~\vec\chi\bigr)\qquad\forall\vec\chi\in\mathbb{V}.\label{eq:cALEweak1}
\end{align}
The detailed derivation of \eqref{eq:cALEweak1} is given in \ref{app:nc}. 
\end{rem}

\subsection{The discrete ALE mappings }\label{sec:DALE}
We follow the discretizations of the interface and bulk domain in \S\ref{sec:ufdis}, but here consider the fitted mesh approach, meaning that the interface mesh is fitted to the bulk mesh such that \[\sigma\in\{\partial o_j^m:\, 1\leq j \leq J_\Omega^m\}\quad\mbox{for all}\quad \sigma\in\{\sigma_j^m:\,1\leq j \leq J_\Gamma\}.\]
Here we keep the mesh connectivity and topology
unchanged, so that in particular $J^m_\Omega = J_\Omega$ and 
$K^m_\Omega = K_\Omega$ for $m=0,\ldots,M$.
Unlike in the unfitted approach, the bulk mesh $\mathscr{T}^m$ can then be divided into simply the interior and exterior elements $\mathscr{T}_-^m$ and $\mathscr{T}_+^m$ with $\mathscr{T}_\Gamma^m = \emptyset$. Therefore, the discrete densities and viscosities can be defined naturally as 
\begin{equation}
\rho^m = \rho_-\mX_{_{\Omega_-^m}}+\rho_+\mX_{_{\Omega_+^m}},\qquad \mu^m = \mu_-\mX_{_{\Omega_-^m}}+\mu_+\mX_{_{\Omega_+^m}},\label{eq:DisQ}
\end{equation}
where $\Omega_-^m$ and $\Omega_+^m$ denote the interior and exterior of $\Gamma^m$, respectively.

For each $m\geq 1$, now assume that we are given the polyhedral surface 
$\Gamma^{m}=\vec X^m(\Gamma^{m-1})$.
We then construct $\mathscr{T}^{m}$ based on $\mathscr{T}^{m-1}$. 
In particular, we update the vertices of the mesh according to
\begin{equation}
\vec a_k^{m} = \vec a_k^{m-1} + \vec\psi^m(\vec a_k^{m-1}),\qquad 1\leq k\leq K_\Omega, \quad 1\leq m\leq M, \label{eq:bdd}
\end{equation}
where $Q^m=\{\vec a_k^m\}_{k=1}^{K_\Omega}$ are the vertices of $\mathscr{T}^m$ and  $\vec\psi^m\in [S_1^{m-1}]^d$ is the displacement of the bulk mesh.  In particular, on introducing 
\begin{align}
\mathbb{Y}^{m-1}&=\bigl\{\vec\chi\in[S_1^{m-1}]^d:\,\vec\chi\cdot\vec n = 0\;\;\mbox{on}\;\;\partial\Omega;\;\vec\chi = \vec X^m - \vec\id\;\;\mbox{on}\;\;\Gamma^{m-1}\bigr\},\nn\\
\mathbb{Y}_0^{m-1}&=\bigl\{\vec\chi\in[S_1^{m-1}]^d:\,\vec\chi\cdot\vec n = 0\;\;\mbox{on}\;\;\partial\Omega;\;\vec\chi = \vec 0\;\;\mbox{on}\;\;\Gamma^{m-1}\bigr\},\nn
\end{align}
we then find $\vec\psi^m\in\mathbb{Y}^{m-1}$ such that 
\begin{equation}
2\bigl(\lambda^{m-1}\,\mat{\tD}(\vec\psi^m),~\mat{\tD}(\vec\chi)\bigr) + \bigl(\lambda^{m-1}\,\nabla\cdot\vec\psi^m,~\nabla\cdot\vec\chi\bigr)=0\qquad\forall\vec\chi\in\mathbb{Y}_0^{m-1},\label{eq:elastic}
\end{equation}
where $\lambda^{m-1} \in S^{m-1}_0$ is defined as
\begin{equation}
\lambda^{m-1}_{|o_j^{m-1}} = 1 + \frac{\max\limits_{o\in\mathscr{T}^{m-1}} |o| - \min\limits_{o\in\mathscr{T}^{m-1}}|o|}{|o_j^{m-1}|},
\qquad j=1,\ldots,J_\Omega.\nn
\end{equation}
The above procedure is used to limit the distortion of small elements 
\cite{Masud1997space,Zhao2020energy}.

Having obtained $\mathscr{T}^{m}$, the discrete mesh velocity 
$\vec W^m \in [S^m_1]^d$ is given by
\begin{equation}
\vec W^m := \sum_{k=1}^{K_\Omega} \left(\frac{\vec a_k^{m} - \vec a_k^{m-1}}{\ttau}\right)\,\phi_k^m,
\label{eq:DMeshV}
\end{equation} 
where $\phi_k^m$ is the nodal basis function of $S_1^m$ at $\vec a_k^m$. The corresponding discrete ALE mappings $\vec{\mathcal{A}}^m[t] \in [S^m_1]^d$, for 
$t\in [t_{m-1},~t_{m}]$, are defined by 
\begin{align}
\vec{\mathcal{A}}^m[t]&: = \vec\id - (t_m - t)\vec W^m=\sum_{k = 1}^{K_\Omega}\left(\frac{t_{m} -t}{\ttau}\,\vec a_k^{m-1}+ \frac{t-t_{m-1}}{\ttau}\,\vec a_k^{m}\right)\phi_k^m.
\label{eq:DALE}
\end{align}
Clearly, $\vec{\mathcal{A}}^m[t_m]$ is the identity map
and the map $\vec{\mathcal{A}}^m[t_{m-1}]\in [S^m_1]^d$ satisfies
\begin{equation} 
\vec{\mathcal{A}}^m[t_{m-1}] = \vec\id - \ttau\,\vec W^m\quad\mbox{with}\quad \vec{\mathcal{A}}^m[t_{m-1}](\Omega_\pm^m) = \Omega_\pm^{m-1}\quad\text{and}\quad
\vec{\mathcal{A}}^m[t_{m-1}](\vec a^m_k) = \vec a^{m-1}_k, \quad k =
1,\ldots,K_\Omega.\nn
\end{equation}
We introduce the Jacobian determinant $\mathcal J^m \in S^m_0$ of the element-wise 
linear map $\vec{\mathcal{A}}^m[t_{m-1}]$ as
\begin{equation}
\mathcal J^{m} := {\rm det}(\nabla \vec{\mathcal{A}}^m[t_{m-1}] ) =
{\rm det}(\mat\Id - \ttau\nabla\vec W^m)= 1 - \ttau\,\nabla\cdot\vec W^m + O(\ttau^2).\label{eq:Jacobian}
\end{equation}
Then we have the following lemma
\begin{lem}\label{lem:energy}
Let $\varphi\in L^2(\Omega)$. Then it holds that 
\begin{equation}
\int_{\Omega_\pm^m}\varphi\circ\vec{\mathcal{A}}^m[t_{m-1}]\mathcal J^{m}\,\dL^d 
= \int_{\Omega_\pm^{m-1}}\varphi\,\dL^d.
\end{equation}
\end{lem}
\begin{proof}
The desired result follows directly from 
$\vec{\mathcal{A}}^m[t_{m-1}](\Omega_\pm^m) = \Omega_\pm^{m-1}$, the definition of 
$\mathcal J^{m}$ and the change-of-variables formula.
\end{proof}

\subsection{The ALE structure-preserving method}
With the fitted finite element approximations, we can then propose the following ALE structure-preserving method based on the weak formulation \eqref{eqn:ALEweak}. 
Let $\vec U^0\in\mathbb{U}^0$ be an approximation of the initial fluid 
velocity $\vec u_0$. 
Moreover, let $\Gamma^0$ be a polyhedral approximation of the initial fluid 
interface $\Gamma_0$ and set $\vec X^0 = \vec\id_{| \Gamma^0}$ with $\mathscr{T}^0$ being a regular fitted partition of $\Omega$. We also set $\Gamma^{-1} = \Gamma^0$ with $\Omega_\pm^{-1}=\Omega_\pm^0$
and $\vec W^0 = \vec 0$, $\mathcal J^0(\vec x) = 1$. Then for $m \geq 0$, find $U^{m+1}\in\mathbb{U}^m$, $P^{m+1}\in\mathbb{P}^m$, $\vec X^{m+1}\in [V^h(\Gamma^m)]^d$ and $\kappa^{m+1}\in V^h(\Gamma^m)$ such that
\begin{subequations}\label{eqn:nALEfem}
\begin{align}
&\bigl(\rho^{m}\,\frac{\vec U^{m+1} - \vec U^m\circ\vec{\mathcal{A}}^m[t_{m-1}]\sqrt{\mathcal J^{m}}}{\ttau},~\vec\chi^h\bigr)+\mathscr{A}(\rho^m, \vec U^m\circ\vec{\mathcal{A}}^m[t_{m-1}]-\vec W^m;~\vec U^{m+1},\vec\chi^h\big)\nn\\&\hspace{1cm}+2\bigl(\mu^{m}\mat{\tD}(\vec U^{m+1}),~\mat{\tD}(\vec\chi^h)\bigr) - \bigl(P^{m+1},~\nabla\cdot\vec\chi^h\bigr)- \gamma\big\langle\kappa^{m+1}\,\vec\nu^m,~\vec\chi^h\big\rangle_{\Gamma^m}\nn\\&\hspace{3cm}= \bigl(\rho^m\vec g,~\vec\chi^h\bigr)\qquad\forall\vec\chi^h\in\mathbb{U}^m,\label{eq:nALEfem1}\\[0.5em]
&\hspace{0.5cm}\bigl(\nabla\cdot\vec U^{m+1},~q^h\bigr) = 0\qquad\forall q^h \in \mathbb{P}^m, \label{eq:nALEfem2}\\[0.5em]
&\hspace{0.5cm}\big\langle\frac{\vec X^{m+1}-\vec\id}{\ttau}\cdot\vec\nu^{m+\frac{1}{2}},~\varphi^h\big\rangle_{\Gamma^m}^h - \big\langle\vec U^{m+1}\cdot\vec\nu^m,~\varphi^h\big\rangle_{\Gamma^m} =0\qquad\forall\varphi^h\in V^h(\Gamma^m),\label{eq:nALEfem3}\\[0.5em]
&\hspace{0.5cm}\big\langle\kappa^{m+1}\,\vec\nu^{m+\frac{1}{2}},~\vec\zeta^h\big\rangle_{\Gamma^m}^h + \big\langle\nabla_s\vec X^{m+1},~\nabla_s\vec\zeta^h\big\rangle_{\Gamma^m} =0\qquad\forall\vec\zeta^h\in [V^h(\Gamma^m)]^d,\label{eq:nALEfem4}
\end{align}
\end{subequations}
and we then set $\Gamma^{m+1} := \vec X^{m+1}(\Gamma^m)$ to construct the new bulk mesh $\mathscr{T}^{m+1}$ through \eqref{eq:bdd}--\eqref{eq:elastic}, and compute the new mesh velocity $\vec W^{m+1}$ through \eqref{eq:DMeshV}. Here on recalling \eqref{eq:Jacobian}, we note that the first term in \eqref{eq:nALEfem1} can be rewritten as
\begin{align}
&\bigl(\rho^m\,\frac{\vec U^{m+1} - \vec U^{m}\circ\vec{\mathcal{A}}^m[t_{m-1}]\sqrt{\mathcal J^{m}}}{\ttau},~\vec\chi^h\bigr)\nn\\&\quad
 = \bigl(\rho^m\frac{\vec U^{m+1} - \vec U^{m}\circ\vec{\mathcal{A}}^m[t_{m-1}]}{\ttau},~\vec\chi^h\bigr) +\tfrac{1}{2} \bigl(\rho^m\,\nabla\cdot\vec W^m,~(\vec U^{m}\circ\vec{\mathcal{A}}^m[t_{m-1}])\cdot\vec\chi^h\bigr) + O(\ttau),\nn
\end{align}
which is hence a consistent temporal discretization of the first two terms in \eqref{eq:ALEweak1}. We have the following theorem for the introduced method \eqref{eqn:nALEfem}, which mimics \eqref{eq:VolumeCon} and \eqref{eq:EnergyLaw} on the discrete level. 
\begin{thm}[stability and volume conservation]
Let $(\vec U^{m+1},P^{m+1},\vec X^{m+1}, \kappa^{m+1})$ be a solution of \eqref{eqn:nALEfem}. Then it holds 
\begin{align}
\mathcal{E}(\rho^k,\vec U^{k+1},\Gamma^{k+1}) + 2\ttau\,\sum_{m=0}^k\norm{\sqrt{\mu^m}\mat{\tD}(\vec U^{m+1})} \leq \mathcal{E}(\rho^{0},\vec U^0,\Gamma^0) + \ttau\,\sum_{m=0}^k\bigl(\rho^{m}\,\vec g,~\vec U^{m+1}\bigr),\label{eq:ALEgDec}
\end{align}
for $k=0,1,\ldots,M-1$.
Moreover, it holds for $m = 0,\ldots, M-1$ that if
\begin{equation}
\left(\mX_{_{\Omega_-^m}}-\omega^m\right)\in\mathbb{P}^m\qquad{\rm with}\quad \omega^m =\frac{\vol(\Omega_-^m)}{\vol(\Omega)},\label{eq:ALEAsup2}
\end{equation}
then 
\begin{equation}
\vol(\Omega_-^{m+1}) = \vol(\Omega_-^m).
\label{eq:DvcALE}
\end{equation}
\end{thm}
\begin{proof}
The proof of \eqref{eq:DvcALE} is exactly the same as that of \eqref{eq:Dvc}
in Theorem~\ref{thm:stavc}. 

For the stability proof, we choose $\vec\chi^h = \ttau\,\vec U^{m+1}$ in \eqref{eq:nALEfem1}, $q^h = P^{m+1}$ in \eqref{eq:nALEfem2}, $\varphi^h = \ttau\,\gamma\kappa^{m+1}$ in \eqref{eq:nALEfem3} and $\vec\zeta^h = (\vec X^{m+1} - \vec\id_{|\Gamma^m})$ in \eqref{eq:nALEfem4}. Combining these equations then gives rise to
\begin{align}
&\bigl(\rho^m\,\delta\vec U^{m},~\vec U^{m+1}\bigr) + 2\ttau\bigl(\mu^m\,\mat{\tD}(\vec U^{m+1}),~\mat{\tD}(\vec U^{m+1})\bigr) 
+ \gamma\big\langle\nabla_s\vec X^{m+1},~\nabla_s(\vec X^{m+1}-\vec\id)\big\rangle_{\Gamma^m} = \ttau\bigl(\rho^m\,\vec g,~\vec U^{m+1}\bigr),\label{eq:nVE}
\end{align}
where \[\delta\vec U^m=\vec U^{m+1}-\vec U^m\circ\vec{\mathcal{A}}^m[t_{m-1}]\sqrt{\mathcal J^m}.\]

Let $(\cdot, \cdot)_{\Omega_\pm^m}$ be the $L^2$-inner product over $\Omega_\pm^m$, respectively.  By the identity $2\vec a\cdot(\vec a - \vec b) = |\vec a|^2 - |\vec b|^2 + |\vec a - \vec b|^2$, we can compute 
\begin{align}
\bigl(\delta\vec U^m,~\vec U^{m+1}\bigr)_{\Omega_\pm^m} &\geq \tfrac{1}{2}\bigl(|\vec U^{m+1}|^2,~1\bigr)_{\Omega_\pm^m} - \tfrac{1}{2}\bigl(|\vec U^{m}\circ\vec{\mathcal{A}}^m[t_{m-1}]|^2,~\mathcal J^{m}\bigr)_{\Omega_\pm^m} \nn\\
&= \tfrac{1}{2}\bigl(|\vec U^{m+1}|^2,~1\bigr)_{\Omega_\pm^m} - \tfrac{1}{2}\bigl(|\vec U^{m}|^2,~1\bigr)_{\Omega_\pm^{m-1}},\label{eq:nVE1}
\end{align}
where the last equality follows from Lemma \ref{lem:energy}. Combing \eqref{eq:nVE1} and \eqref{eq:DisQ} yields
\begin{align}
&\bigl(\rho^m\,\delta\vec U^m,~\vec U^{m+1}\bigr) = \rho_-\bigl(\delta \vec U^m,~\vec U^{m+1}\bigr)_{\Omega_-^m} +\rho_+\bigl(\delta \vec U^m,~\vec U^{m+1}\bigr)_{\Omega_+^m}\nn\\ & \quad
\geq \tfrac12 \rho_-\bigl(|\vec U^{m+1}|^2,~1\bigr)_{\Omega_-^m} - \tfrac12 \rho_-(|\vec U^{m}|^2,~1\bigr)_{\Omega_-^{m-1}}+\tfrac12 \rho_+\bigl(|\vec U^{m+1}|^2,~1\bigr)_{\Omega_+^m} - \tfrac12 \rho_+(|\vec U^{m}|^2,~1\bigr)_{\Omega_+^{m-1}}\nn\\&\quad
= \tfrac{1}{2}\bigl(\rho^m,~|\vec U^{m+1}|^2\bigr) - \tfrac{1}{2}\bigl(\rho^{m-1},~|\vec U^{m}|^2\bigr).\label{eq:nVE3}
\end{align}
Now inserting \eqref{eq:nVE3} into \eqref{eq:nVE}, and recalling \eqref{eq:Elem}, we obtain 
\begin{equation}
\mathcal{E}(\rho^m,\vec U^{m+1}, \Gamma^{m+1}) + 2\ttau\norm{\sqrt{\mu^m}\,\mat{\tD}(\vec U^{m+1})} \leq \mathcal{E}(\rho^{m-1}, \vec U^m, \Gamma^m) + \ttau\bigl(\rho^m\vec g,~\vec U^{m+1}\bigr).
\label{eq:ALEDes}
\end{equation}
Summing \eqref{eq:ALEDes} for $m=0,\ldots,k$ yields \eqref{eq:ALEgDec},
on recalling that $\Omega_\pm^{-1}=\Omega_\pm^0$.
\end{proof}
 
\begin{rem} The condition \eqref{eq:ALEAsup2} trivially holds if the pair elements \eqref{eq:P2P0} or \eqref{eq:P2P1P0} are used for the fluid flow in the bulk. This then leads to the exact volume preservation property. In practice, the moving mesh approach via \eqref{eq:elastic} in general works smoothly (see {\bf Example 4} in \S \ref{sec:num}). Nevertheless, in the case of very large deformations (see {\bf Example 5} in \S\ref{sec:num}), a remeshing of the bulk mesh may become necessary. Then the obtained velocity solution needs to be appropriately interpolated onto the new mesh and the stability result \eqref{eq:ALEgDec} will in general no longer hold. Besides, we note that \eqref{eq:nALEfem4} will lead to a very good property of the interface mesh \cite{Barrett20}, thus no remeshing for the interface mesh is necessary in general.  
\end{rem}
\begin{rem} 
We note that \eqref{eqn:nALEfem} is similar to \cite[(2.26)--(2.29)]{Duan2022energy} except that here we employ the semi-implicit approximation
of the unit normal from \eqref{eq:weightv} instead of $\vec\nu^m$ for the stability and volume conservation. In addition, the scheme from \cite{Duan2022energy} is based on a formulation with a divergence free ALE mesh velocity, see Remark \ref{rem:DL}. In the construction of the discrete ALE mappings, the discrete displacement is then solved via a Stokes equation, instead of the elastic equation we consider in \eqref{eq:elastic}.  
\end{rem}
 
\begin{rem}On recalling Remark \ref{rem:cALE}, it is also possible to consider a discretization of the conservative ALE formulation \eqref{eq:cALEweak1}. We introduce 
\begin{equation}
\Omega_\pm^h(t) = \vec{\mathcal{A}}^m[t](\Omega_\pm^m),\quad t\in[t_{m-1}, t_{m}].\nn
\end{equation}
It is easy to show that (see \cite{Elliott2012ale})
\begin{equation}
\partial_t^\circ(\vec\chi^h\circ\vec{\mathcal{A}}^m[t]^{-1}) = \vec 0\quad\mbox{for all}\quad\vec\chi^h\in\mathbb{U}^m.\nn
\end{equation}
Then {\bf the conservative ALE method} based on \eqref{eq:cALEweak1} is 
the same as \eqref{eqn:nALEfem}, with  \eqref{eq:nALEfem1} replaced by
\begin{align}
&\frac{1}{\ttau}\bigl[\bigl(\rho^{m}\,\vec U^{m+1},~\vec\chi^h\bigr) -\bigl(\rho^{m-1}\,\vec U^{m},~\vec\chi^h\circ\vec{\mathcal{A}}^m[t_{m-1}]^{-1}\bigr)\bigr]-\frac{1}{2\,\ttau}\mathscr{B}^m(\rho^m,\vec W^m,~\vec U^{m+1}, \vec\chi^h) \nn\\
&\hspace{1cm} +\mathscr{A}(\rho^m, \vec U^{m}\circ\vec{\mathcal{A}}^m[t_{m-1}]-\vec W^m;~\vec U^{m+1},\vec\chi^h\big)+2\bigl(\mu^{m}\,\mat{\tD}(\vec U^{m+1}),~\mat{\tD}(\vec\chi^h)\bigr)\nn\\&\hspace{2cm} - \bigl(P^{m+1},~\nabla\cdot\vec\chi^h\bigr) - \gamma\big\langle\kappa^{m+1}\,\vec\nu^m,~\vec\chi^h\big\rangle_{\Gamma^m} = \bigl(\rho^m\vec g,~\vec\chi^h\bigr)\qquad\forall\vec\chi^h\in\mathbb{U}^m,\label{eq:cALEfem1}
\end{align}
where we introduced the time-integrated term (see \cite{Liu2013,ivanvcic2022energy})
\begin{equation}
\mathscr{B}^m(\rho^m,\vec W^m,~\vec U^{m+1}, \vec\chi^h) = \int_{t_{m-1}}^{t_m}\bigl([\rho^m\circ\vec{\mathcal{A}}^m[t]^{-1}]\,~\nabla\cdot[\vec W^m\circ\vec{\mathcal{A}}^m[t]^{-1}],~[\vec U^{m+1}\cdot\vec\chi^h]\circ\vec{\mathcal{A}}^m[t]^{-1}\bigr)\,\rd t.\nn
\end{equation}
The proof of the volume conservation for the above method is straightforward. Using the technique in \cite{ivanvcic2022energy}, it is also not difficult to show the proposed method also satisfies the energy stability in \eqref{eq:ALEgDec}.  
\end{rem}

\section{Numerical tests}
\label{sec:num}
In this section, we present several benchmark tests for the introduced structure-preserving methods \eqref{eqn:fem} and \eqref{eqn:nALEfem} in both 2d and 3d. In what follows, we always choose $\vec U^0 = \vec 0$. 
To solve the weakly nonlinear systems \eqref{eqn:fem}, we use a Picard-type iteration \eqref{eqn:Picardfem} and choose ${\rm tol} = 10^{-8}$ in \eqref{eq:iterativeTol}. The linear systems from \eqref{eqn:Picardfem} can be solved efficiently with the Schur complement approach and preconditioned Krylov iterative solvers  described in \cite{BGN15stable}. For \eqref{eqn:nALEfem}, we employ a similar Picard-type iteration, and the resulting linear systems are solved via a sparse LU factorization with the open library Eigen, see \cite{Eigenweb}. For ease of presentation, we denote by ``\ESP-SP''  the Eulerian structure-preserving method in \eqref{eqn:fem}. We employ a bulk mesh adaptation strategy as described in \cite{BGN15stable} with the same notations ``$n{\rm adapt}_{k,l}$''  to denote $\ttau = 10^{-3}/n$, $N_f= 2^k$ and $N_c=2^l$. For the case $d = 2$
we have that $K_\Gamma=J_\Gamma = 2^k$, while for $d = 3$ it holds that $(K_\Gamma, J_\Gamma)=(770, 1536), (1538, 3072)$ for $k = 5, 6$. Besides, we use ``$\ALE_n$-SP'' to denote the ALE structure-preserving method \eqref{eqn:nALEfem}, and ``$\ALE_c$-SP'' to denote \eqref{eqn:nALEfem} in the case when \eqref{eq:nALEfem1} is replaced by \eqref{eq:cALEfem1}.

\subsection{Convergence tests}

\newcommand{\LerrorPpc}{\|P_c - p_c\|_{L^2}}
\newcommand{\errorLl}{\|\lambda^h - \lambda\|_{L^\infty}}

For the convergence tests, we consider a static bubble \cite{Ganesan07,Gross07extended, BGN2013eliminating} and an expanding bubble \cite{Agnese20}. In particular, we choose  $\Gamma(t)=\{\vecz\in\bR^d: |\vecz|=r(t)\}$ and 
\begin{equation} \label{eq:rupg} 
r(t)=([r(0)]^d + \alpha\,t\,d)^{\frac{1}{d}},\quad \vec u(\vec x, t)= \vec
u_b(\vec x) = 
\alpha\frac{\vec x}{|\vec x|^{d}},\quad p(\vec x, t) = \lambda(t)[\mX_{_{\Omega_-(t)}}-\omega(t)],\quad \vec g = \alpha^2(1-d)\frac{\vec x}{|\vec x|^{2d}},\nn
\end{equation}
where $\alpha\geq 0$, $\omega(t)=\frac{\vol(\Omega_-(t))}{\vol(\Omega)}$ and $\lambda(t)=\gamma(d-1)[r(t)]^{-1}+2\alpha(d-1)(\mu_+-\mu_-)[r(t)]^{-d}$. In the case of $\alpha=0$, this corresponds to the trivial stationary solution of a circle or sphere. When $\alpha>0$, the interface corresponds to an expanding circle or sphere within a domain that does not contain the origin. We also need to replace \eqref{eqn:BD} by imposing the inhomogeneous conditions $\vec u(\cdot, t)=\vec u_b$ on $\partial_1\Omega=\partial\Omega$. Nevertheless, it is not difficult to show that the Eulerian formulation \eqref{eq:weak1} still holds as the test function $\vec\chi=\vec 0$ on $\partial\Omega$. In the ALE approach, on $\partial\Omega$ we impose $\vec w=\vec 0$ instead of \eqref{eq:MeshVBD} so that the domain boundary is fixed. Similarly \eqref{eq:ALEweak1} and \eqref{eq:cALEweak1} still hold since the test function $\vec\chi$ vanishes on the boundary. Moreover, we introduce the errors 
\begin{align}
&\norm{\vec X-\vecz}_{L^\infty}=\max_{m=1,\cdots,M}\norm{\vec X^m-\vecz(\cdot,t_m)}_{L^\infty}\quad\mbox{with}\quad \norm{\vec X^m - \vecz(\cdot,t_m)}_{L^\infty} = \max_{k=1,\cdots, K_\Gamma}\,\min_{\vec y\in\Gamma(t_m)}|\vec X^m(\vec q_k^m)-\vec y|,\nn\\
& \norm{\vec U-I_2^h\vec u}_{L^\infty}=\max_{m=1,\cdots,M}\norm{\vec U^m-I_2^m\vec u(\cdot, t_m)}_{L^\infty},\qquad
\norm{P-p}_{L^2} = \left(\ttau\sum_{m=1}^M\norm{P^m-p(\cdot,t_m)}_{L^2}^2\right)^{\frac{1}{2}},\nn
\end{align}
where for the calculation of $\norm{P^m-p(\cdot,t_m)}_{L^2}$ we use a
quadrature rule that is exact for polynomials of degree 17.
As discretization parameters for the \ESP-SP method we 
choose an adaptive bulk mesh with
$h_c = 8 h_f$ and $h^m_\Gamma \approx h_f$.
For the ALE method, we choose three meshes given by (i) $J_\Gamma=24$, $J_\Omega=184$,  $K_\Omega=110$; (ii) $J_\Gamma=48$, $J_\Omega=684$, $K_\Omega=376$; (iii) $J_\Gamma=96$,  $J_\Omega=2494$, $K_\Omega=1315$.

\begin{table}[!htp]
\centering
\def\temptablewidth{0.79\textwidth}
\vspace{-2pt}
\caption{
 Convergence experiments for the static bubble over the time interval $[0,1]$ with the parameters as in \eqref{eq:stat}. The upper and middle panels are the results from the \ESP-SP method using the pair element P2-P1 without or with XFEM, at the lower panel are the results from the $\ALE_n$-SP method with the pair element P2-(P1+P0), the time step is fixed as $\ttau=10^{-3}$. }\label{tab:stat1}
{\rule{\temptablewidth}{1pt}}
\begin{tabular}{c|cc|cc|cc}
$J_\Gamma$  &$\norm{\vec X-\vecz}_{L^\infty}$ &order & $\norm{\vec U-I_2^h\vec u}_{L^\infty}$  &order & $\norm{P-p}_{L^2}$  &order  \\ \hline  
64  & 2.19E-3 &-    & 4.37E-3 &-    & 2.09E-1 &-    \\ \hline 
128 & 9.72E-4 &1.17 & 2.31E-3 &0.92 & 1.48E-1 &0.50 \\ \hline 
256 & 4.53E-4 &1.10 & 1.19E-3 &0.96 & 1.03E-1 &0.52  
\end{tabular}
{\rule{\temptablewidth}{1pt}}
{\rule{\temptablewidth}{1pt}}
\begin{tabular}{c|cc|cc|cc} 
$J_\Gamma$  &$\norm{\vec X-\vecz}_{L^\infty}$ &order & $\norm{\vec U-I_2^h\vec u}_{L^\infty}$  &order & $\norm{P-p}_{L^2}$ &order  \\ \hline
64  & 0 &- &0 &- & 6.42E-2 &-    \\ \hline 
128 & 0 &- &0 &- & 3.73E-2 &0.78 \\ \hline 
256 & 0 &- &0 &- & 1.62E-2 &1.20  
\end{tabular}
{\rule{\temptablewidth}{1pt}}
{\rule{\temptablewidth}{1pt}}
\begin{tabular}{c|cc|cc|cc} 
$J_\Gamma$  &$\norm{\vec X-\vecz}_{L^\infty}$ &order & $\norm{\vec U-I_2^h\vec u}_{L^\infty}$  &order & $\norm{P-p}_{L^2}$  &order  \\ \hline
24 &0 &- &0 &- &1.90E-1 &- \\ \hline 
48 &0 &- &0 &- &9.08E-2 &1.07 \\ \hline 
96 &0 &- &0 &- &3.15E-2 &1.53
\end{tabular}
{\rule{\temptablewidth}{1pt}}
\end{table}

\begin{figure}[!htp]
\centering
\includegraphics[width=0.45\textwidth]{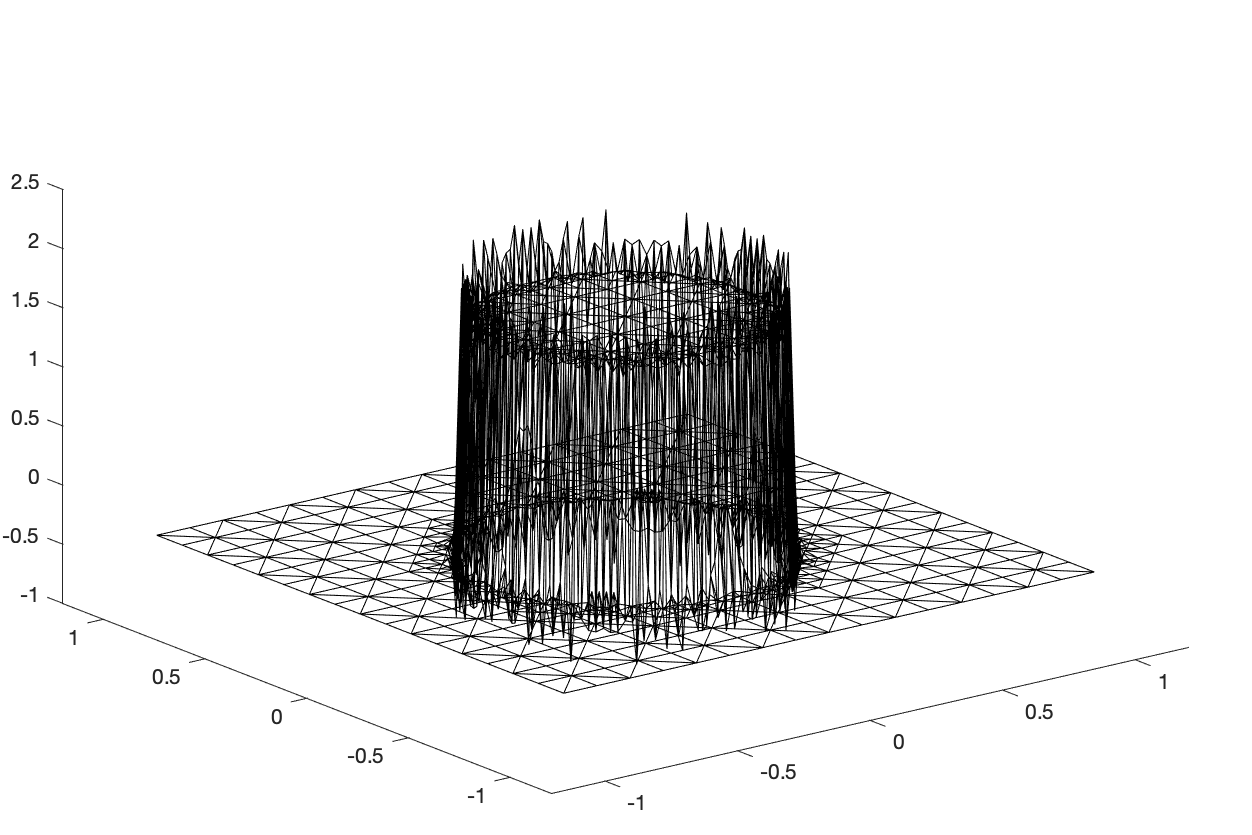}\hspace{0.05cm}
\includegraphics[width=0.45\textwidth]{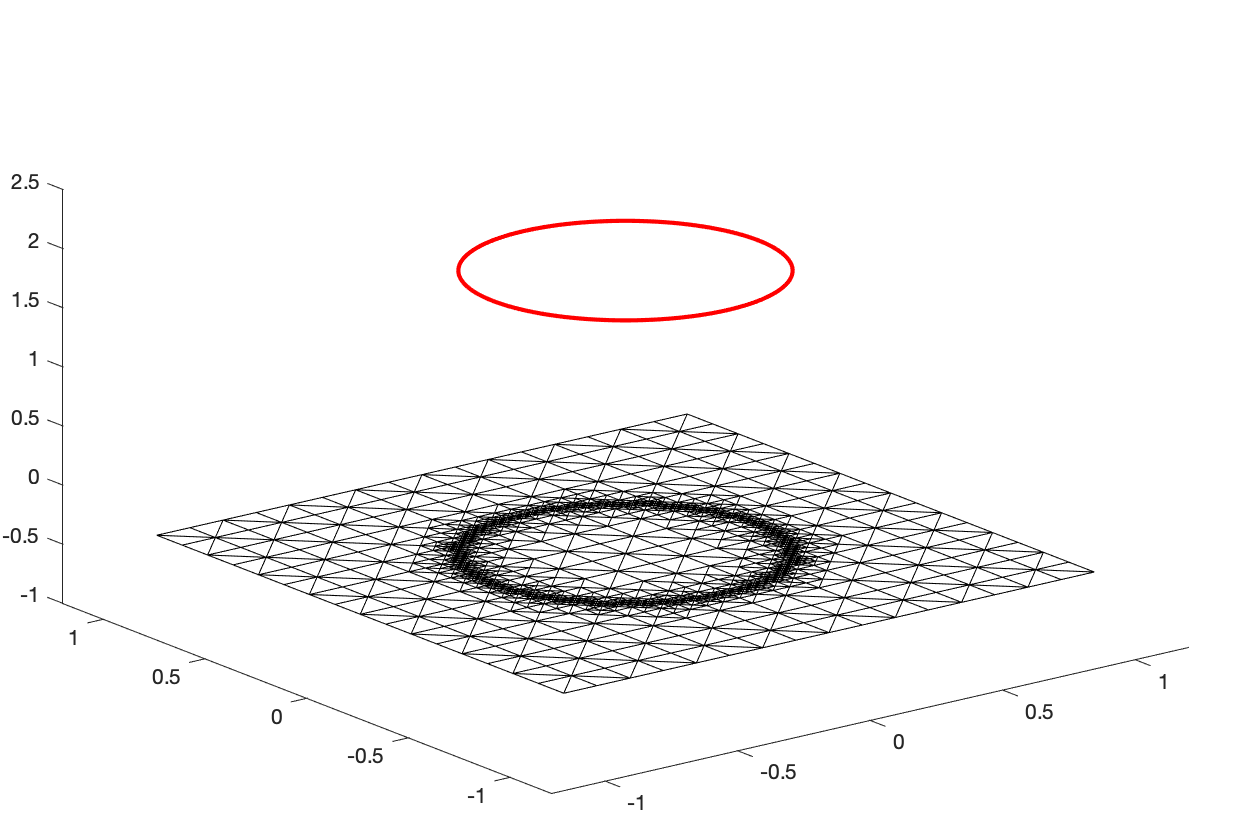}
\includegraphics[width=0.45\textwidth]{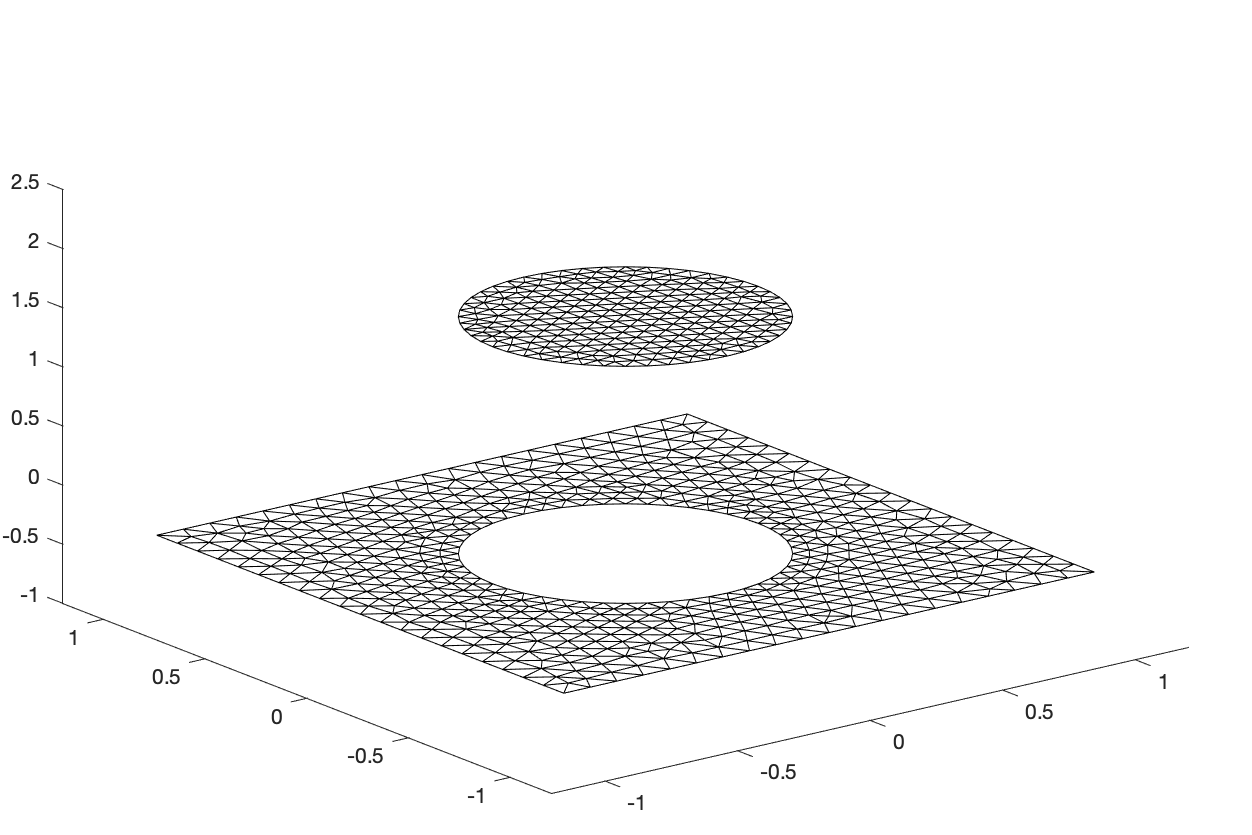}
\caption{Discrete pressure plots for the static bubble at time 
$t=1$, corresponding to the rows with $J_\Gamma=128$ and $J_\Gamma=48$
in Table~\ref{tab:stat1}.
Top left panel: the \ESP-SP method without XFEM, using the spaces 
\eqref{eq:P2P1}. 
Top right panel: the \ESP-SP method with XFEM, using the spaces 
\eqref{eq:P2P1XFEM}. 
Here we plot the $S^m_1$ part and the $\mX_{\Omega_-^m}$ part of the
discrete pressure separately, so that the pressure in the inner phase is the
sum of the two. 
Lower panel: the $\ALE_n$-SP method using the spaces \eqref{eq:P2P1P0}.}
\label{fig:pres}
\end{figure}

\begin{table}[!htp]
\centering
\def\temptablewidth{0.79\textwidth}
\vspace{-2pt}
\caption{Convergence experiments for the expanding bubble over the time interval $[0,1]$ with parameters as in \eqref{eq:epd}. The upper and middle panels are the results from the \ESP-SP method using the pair element P2-P1 without or  with XFEM, at the lower panel are the results from the $\ALE_n$-SP method with the pair element P2-(P1+P0).  }\label{tab:epd1}
{\rule{\temptablewidth}{1pt}}
\begin{tabular}{cc|cc|cc|cc}
$J_\Gamma$ & $\ttau$  &$\norm{\vec X-\vecz}_{L^\infty}$ &order & $\norm{\vec U-I_2^h\vec u}_{L^\infty}$  &order & $\norm{P-p}_{L^2}$  &order  \\ \hline  
 64 & $10^{-2}$ & 2.43E-4 &-    & 1.28E-2 &-    & 1.45E-0 &-    \\ \hline 
128 & $10^{-3}$ & 1.53E-4 &0.67 & 7.16E-3 &0.84 & 9.71E-1 &0.58 \\ \hline 
256 & $10^{-4}$ & 6.61E-5 &1.21 & 4.67E-3 &0.62 & 6.90E-1 &0.49 
\end{tabular}
{\rule{\temptablewidth}{1pt}}
{\rule{\temptablewidth}{1pt}}
\begin{tabular}{cc|cc|cc|cc} 
$J_\Gamma$ & $\ttau$  &$\norm{\vec X-\vecz}_{L^\infty}$ &order & $\norm{\vec U-I_2^h\vec u}_{L^\infty}$  &order & $\norm{P-p}_{L^2}$ & order \\ \hline
 64 & $10^{-2}$ & 4.72E-4 &-    & 1.28E-2 &-    & 1.32E-0  &-    \\ \hline 
128 & $10^{-3}$ & 1.18E-4 &2.00 & 6.91E-3 &0.89 & 8.74E-1 &0.59 \\ \hline 
256 & $10^{-4}$ & 3.11E-5 &1.92 & 4.45E-3 &0.63 & 6.08E-1 &0.52 
\end{tabular}
{\rule{\temptablewidth}{1pt}}
{\rule{\temptablewidth}{1pt}}
\begin{tabular}{cc|cc|cc|cc} 
$J_\Gamma$ & $\ttau$  &$\norm{\vec X-\vecz}_{L^\infty}$ &order & $\norm{\vec U-I_2^h\vec u}_{L^\infty}$  &order & $\norm{P-p}_{L^2}$ & order \\ \hline
24 & $10^{-2}$ &2.59E-3 &-  &1.08E-2 &-  & 9.02E-1  &- \\ \hline 
48 & $10^{-3}$ &6.57E-4  &1.98 &5.56E-3  &0.96 &4.08E-1  &1.14 \\ \hline 
96 & $10^{-4}$ &1.70E-4  &1.95 & 1.24E-3 &2.16  &1.68E-1  &1.28
\end{tabular}
{\rule{\temptablewidth}{1pt}}
\end{table}

\begin{figure}[!htp]
\centering
\includegraphics[width=0.45\textwidth]{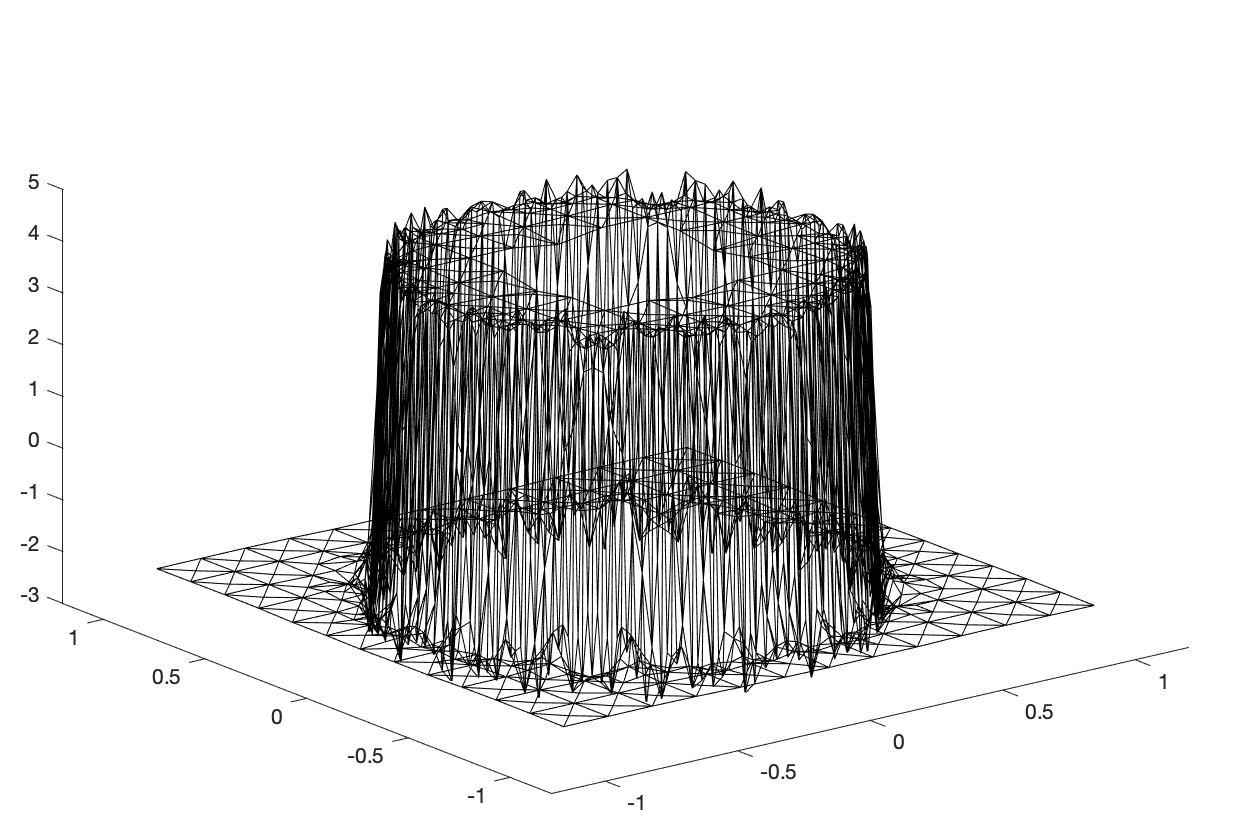}\hspace{0.05cm}
\includegraphics[width=0.45\textwidth]{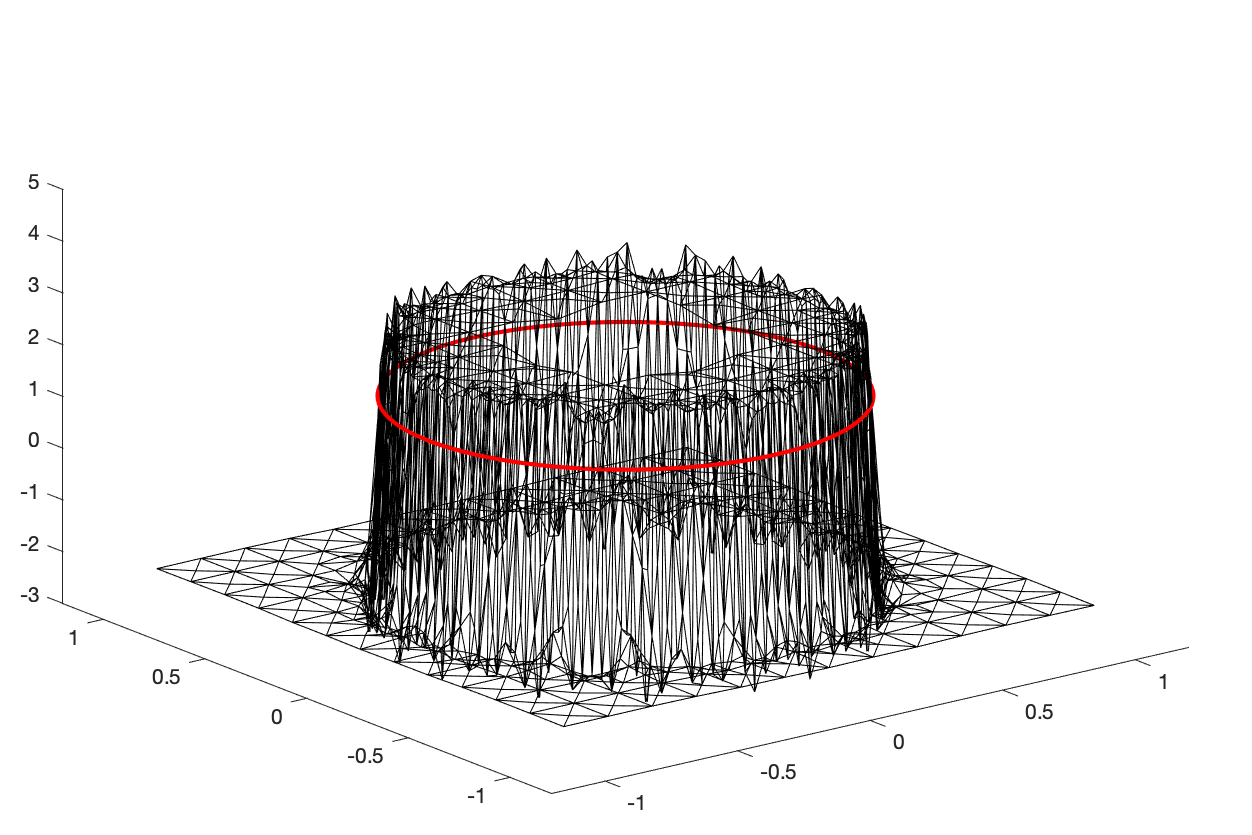}
\includegraphics[width=0.45\textwidth]{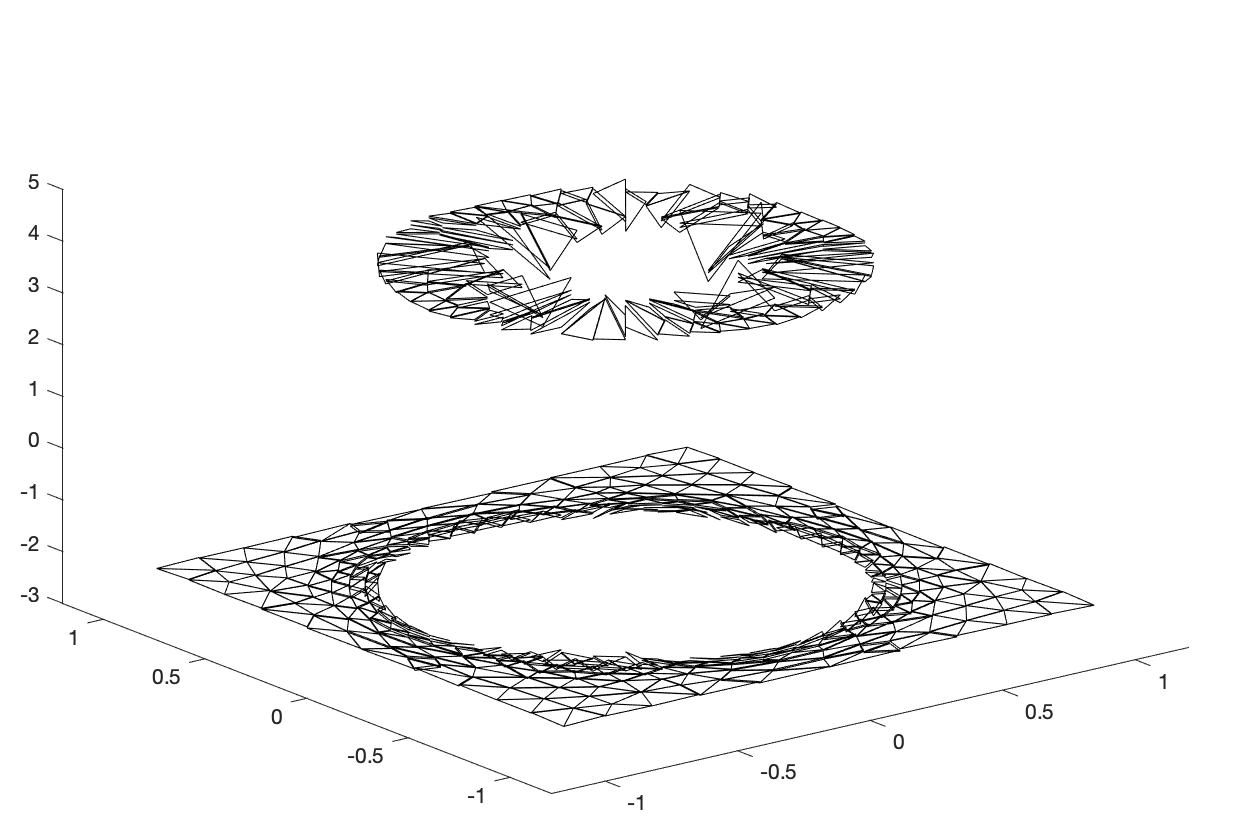}
\caption{Discrete pressure plots for the expanding bubble at time 
$t=1$, corresponding to the rows with $J_\Gamma=128$ and $J_\Gamma=48$
in Table~\ref{tab:epd1}.
Top left panel: the \ESP-SP method without XFEM, using the spaces 
\eqref{eq:P2P1}. 
Top right panel: the \ESP-SP method with XFEM, using the spaces 
\eqref{eq:P2P1XFEM}. 
Here we plot the $S^m_1$ part and the $\mX_{\Omega_-^m}$ part of the
discrete pressure separately, so that the pressure in the inner phase is the
sum of the two. 
Lower panel: the $\ALE_n$-SP method using the spaces \eqref{eq:P2P1P0}.}
\label{fig:pres2}
\end{figure}

\vspace{0.6em}\noindent
{\bf Example 1}:  For the stationary solution, we consider $\Omega=(-1,1)^2$ and choose 
\begin{equation}
\alpha=0,\quad\gamma = 1,\quad \rho_+=1000,\quad\rho_-=100, \quad\mu_\pm = 1,\quad r(0)=0.5.\label{eq:stat}
\end{equation}
The numerical results obtained by the \ESP-SP method and the $\ALE_n$-SP method are reported in Table \ref{tab:stat1}. Here we observe that the enriched \ESP-SP
method and the $\ALE_n$-SP method manage to capture the zero bulk velocity
exactly, see also \cite{BGN2013eliminating,Agnese16} for more details. Moreover, the jump of the pressure across the interface can be accurately captured by the enriched \ESP-SP method and the $\ALE_n$ method, as shown in Fig.~\ref{fig:pres}. However, for the \ESP-SP method without XFEM, we observe that the jump of the pressure is not well resolved by the continuous piecewise linear P1 element. This results in a decrease of the convergence rate for the pressure errors, as illustrated in Table \ref{tab:stat1}.

\vspace{0.6em}\noindent
{\bf Example 2}:  For the expanding bubble, we consider $ \Omega=(-1,1)^2\setminus[-\frac{1}{3}, \frac{1}{3}]^2$ and choose 
\begin{equation}
\alpha=0.15,\quad\gamma=1, \quad \rho_\pm=100,\quad\mu_+=10,\quad\mu_-=1,\quad r(0)=0.5.\label{eq:epd}
\end{equation}
 The numerical results are reported in Table \ref{tab:epd1}. For the \ESP-SP method without XFEM, we still observe that the order of convergence for the pressure errors is about 0.5. However, similar pressure errors are observed for the enriched \ESP-SP method. This is likely due to the inaccurate capture of the jump in the viscosities, which introduces additional errors in the case of a nonzero velocity. Nevertheless, a second order convergence is still observed for the interface errors in the enriched \ESP-SP method, similar to the $\ALE_n$-SP method. The discrete pressure plots for the expanding bubble are shown in Fig.~\ref{fig:pres2}.

\subsection{The rising bubble}

We also study the dynamics of a rising bubble in two different cases,  which were considered in \cite{Hysing2009}. The physical parameters are given by 
\begin{itemize}
\item {\bf Case I}:
\begin{equation}
\rho_+=1000,\quad \rho_- = 100,\quad \mu_+=10, \quad \mu_-= 1,\quad \gamma = 24.5,\quad\vec g = -0.98\vec e_d;\label{eq:Bench1}
\end{equation}
\item {\bf Case II}:
\begin{equation}
\rho_+=1000,\quad \rho_- = 1,\quad \mu_+=10, \quad \mu_-= 0.1,\quad \gamma = 1.96,\quad\vec g = -0.98\vec e_d,\label{eq:Bench2}
\end{equation} 
\end{itemize}
where $\vec e_d = (0,1)^T$ in 2d and $(0,0,1)^T$ in 3d.
We define the following discrete benchmark quantities 
\begin{alignat}{2}
\cira|_{t_m} &:=\frac{\pi^{\frac{1}{d}}[2d\vol(\Omega_-^m)]^{\frac{d-1}{d}}}{|\Gamma^m|} ,\qquad V_c|_{t_m}&&:=\frac{\int_{\Omega_-^m}(\vec U^m\cdot\vec e_d)\,\dL^d}{\vol(\Omega_-^m)},\nn\\
y_c|_{t_m}&:=\frac{\int_{\Omega_-^m}(\vec\id\cdot\vec e_d)\,\dL^d}{\vol(\Omega_-^m)}, \qquad\quad   v_\Delta|_{t_m}&&: = \frac{\vol(\Omega_-^m)- \vol(\Omega_-^0)}{\vol(\Omega_-^0)},\nn
\end{alignat}
where $\cira$ denotes the degree of circularity (sphericity) of the bubble, $V_c$ is the bubble's rise velocity, $y_c$ is the bubble's centre of mass in the vertical direction, and $v_\Delta$ is the relative volume loss.

\subsubsection{Numerical results in 2d} \label{sec:nr2d}

We consider the problem in a bounded domain $\overline{\Omega} = [0, 1]\times[0,2]$ with $\partial_1\Omega=[0,1]\times\{0,2\}$ and $\partial_2\Omega=\{0,1\}\times[0,2]$. Moreover, the initial interface is given by $\Gamma(0):=\bigl\{\vecz\in\Omega:\;|\vecz - (\frac{1}{2},\frac{1}{2})^T|=\frac{1}{4}\bigr\}$. 

\begin{table}[!htp]
\centering
\def\temptablewidth{0.7\textwidth}
\vspace{-4pt}
\caption{2d benchmark quantities of the rising bubble in case I from the Eulerian-SP method, where we use the pair element P2-P1 with XFEM. Here ``${\rm Hysing3}$'' denotes the finest discretization run of group 3 in \cite{Hysing2009}. }\label{tab:Euler1}
{\rule{\temptablewidth}{1pt}}
\begin{tabular}{c|cccc|cc}
  &${\rm adapt}_{5,2}$  & ${\rm adapt}_{7,3}$ & 2${\rm adapt}_{9,4}$ &5${\rm adapt}_{11,5}$  &${\rm Hysing3}$\\\hline  
$\cira_{\min}$ & 0.9135 &0.9068 & 0.9034  &0.9019 &0.9013\\\hline 
$t_{\cira = \cira_{\min}}$ &2.0770 &1.9420 & 1.9110 &1.9028  &1.9000 \\[0.4em]\hline 
$V_{c,\max}$ &0.2479 &0.2415 &0.2414 &0.2416 &0.2417\\\hline 
$t_{_{V_c = V_{c,\max}}}$ &0.9470 &0.9360 & 0.9255 &0.9200 &0.9239\\\hline 
$y_c(t=3)$ &1.0907 &  1.0823 & 1.0815 &1.0817 &1.0817 
\end{tabular}
{\rule{\temptablewidth}{1pt}}
\end{table}

\begin{figure}[!htp]
\centering
\includegraphics[width=0.9\textwidth]
{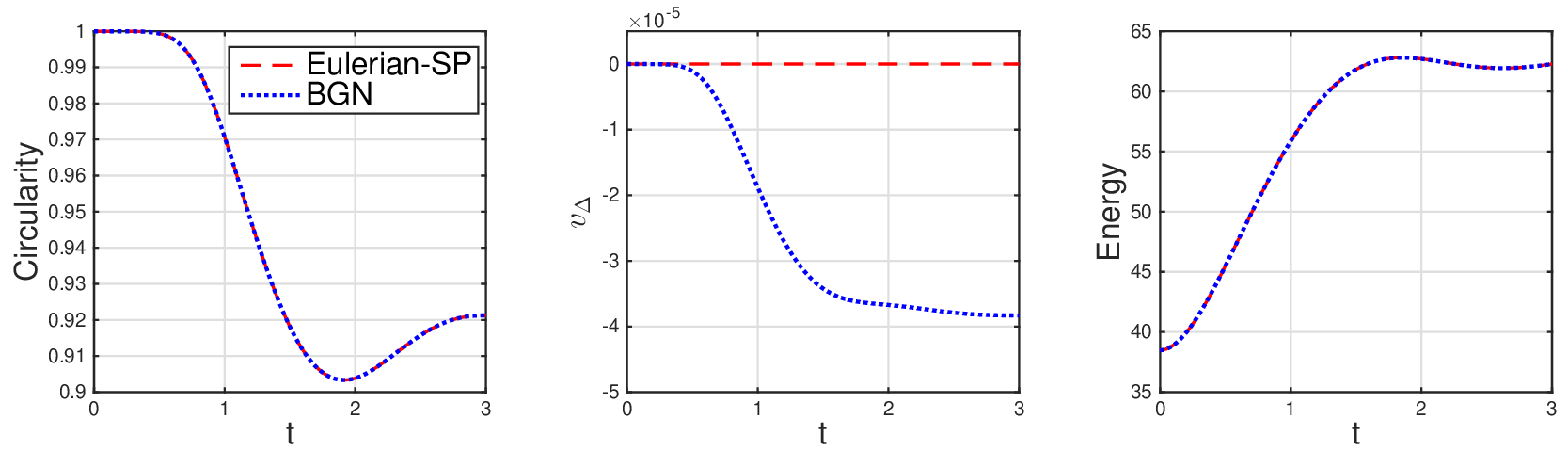}
\caption{ Plots of the circularity of the rising bubble, the relative volume loss and the discrete energy over time for case I in 2d,
using $2{\rm  adapt}_{9,4}$, where ``BGN'' refers to the linear stable scheme introduced in \cite{BGN15stable}.}
\label{fig:EulerBGNC}
\end{figure}

\vspace{0.6em}\noindent
{\bf Example 3}: We start with the rising bubble of case I and consider the Eulerian-SP method \eqref{eqn:fem}. We employ the finite element spaces P2-P1 in \eqref{eq:P2P1} with XFEM so that \eqref{eq:Asup2} is satisfied, see Remark \ref{rem:ufXFEM}.  The quantitative benchmark quantities for the rising bubble are reported in Table \ref{tab:Euler1}. As a comparison, we include the results from the finest discretization run of group 3 in \cite{Hysing2009} (denoted by ``Hysing3''), which can provide very accurate and consistent approximations based on the various benchmark tests e.g., \cite{Aland2012,BGN15stable,Frachon2019cut, Agnese20, Duan2022energy}.  We observe that these quantitative values are in good agreement with corresponding results of ``Hysing3''. 

We next compare the introduced \ESP-SP method with the BGN method in \cite{BGN15stable}. We find that the results in Table \ref{tab:Euler1} are quite similar to those in \cite[Table 2]{BGN15stable}. Besides, we plot the time history of the circularity $\cira$, the relative volume loss $v_\Delta$ and the discrete energy for the two methods in Fig.~\ref{fig:EulerBGNC}. We observe that the evolution of the circularity and energy shows very good agreement between the two methods. Nevertheless, the introduced \ESP-SP method can exactly preserve the enclosed volume of the bubble, while the BGN method does not. This numerically confirms the volume preserving property of the \ESP-SP method, see \eqref{eq:Dvc}.  

\begin{table}[!htp]
\centering
\def\temptablewidth{0.8\textwidth}
\vspace{-4pt}
\caption{2d benchmark quantities of the rising bubble in case I from the ALE structure-preserving methods. Here we denote $h=1/J_\Gamma$ and use the pair element P2-P0 (upper panel) and P2-(P1+P0) (lower panel), $h_0 = 1/32$ and $\ttau_0=0.01$.}
\label{tab:ALEBench}
{\rule{\temptablewidth}{1pt}}
\begin{tabular}{c|ccc|ccc}
& \multicolumn{3}{c|}{${\rm ALE}_n$-SP} & \multicolumn{3}{c}{${\rm ALE}_c$-SP}  \\\hline
 ($h$,~$\ttau$) & $(h_0, \ttau_0)$ & $(\frac{h_0}{2}, \frac{\ttau_0}{4})$ & $(\frac{h_0}{4}, \frac{\ttau_0}{16})$
 & $(h_0, \ttau_0)$ & $(\frac{h_0}{2}, \frac{\ttau_0}{4})$ & $(\frac{h_0}{4}, \frac{\ttau_0}{16})$  \\ \hline
$\cira_{\min}$   & 0.9221 & 0.9084 & 0.9032  & 0.9222 & 0.9084 & 0.9032 \\ \hline 
$t_{\cira = \cira_{\min}}$ & 1.8400 & 1.8800 & 1.9019  & 1.8500 & 1.8825 & 1.9025 \\ \hline 
$V_{c,\max}$                     & 0.2246 & 0.2355 & 0.2396  & 0.2255 & 0.2356 & 0.2396\\ \hline 
$t_{V_c = V_{c,\max}}$           & 0.9800 & 0.9450 & 0.9306  & 0.9800 & 0.9425 & 0.9306 \\ \hline 
$y_c(t=3)$                       & 1.0855 & 1.0811 & 1.0805  & 1.0835 & 1.0806 & 1.0803
\end{tabular}
{\rule{\temptablewidth}{1pt}}
{\rule{\temptablewidth}{1pt}}
\begin{tabular}{c|ccc|ccc}
 ($h$,~$\ttau$) & $(h_0, \ttau_0)$ & $(\frac{h_0}{2}, \frac{\ttau_0}{4})$ & $(\frac{h_0}{4}, \frac{\ttau_0}{16})$
 & $(h_0, \ttau_0)$ & $(\frac{h_0}{2}, \frac{\ttau_0}{4})$ & $(\frac{h_0}{4}, \frac{\ttau_0}{16})$  \\ \hline
$\cira_{\min}$   & 0.9030 & 0.9024 & 0.9015  & 0.9031 & 0.9023 & 0.9015 \\ \hline 
$t_{\cira = \cira_{\min}}$ & 1.9400 & 1.9050 & 1.9000  & 1.9500 & 1.9075 & 1.9006 \\ \hline 
$V_{c,\max}$                     & 0.2440 & 0.2423 & 0.2418  & 0.2440 & 0.2423 & 0.2418\\ \hline 
$t_{V_c = V_{c,\max}}$           & 0.9300 & 0.9250 & 0.9225  & 0.9300 & 0.9250 & 0.9225 \\ \hline 
$y_c(t=3)$                       & 1.0905 & 1.0840 & 1.0823  & 1.0885 & 1.0836 & 1.0822 
\end{tabular}
{\rule{\temptablewidth}{1pt}}
\end{table}%

\begin{figure}[!htp]
\centering
\includegraphics[width=0.8\textwidth]
{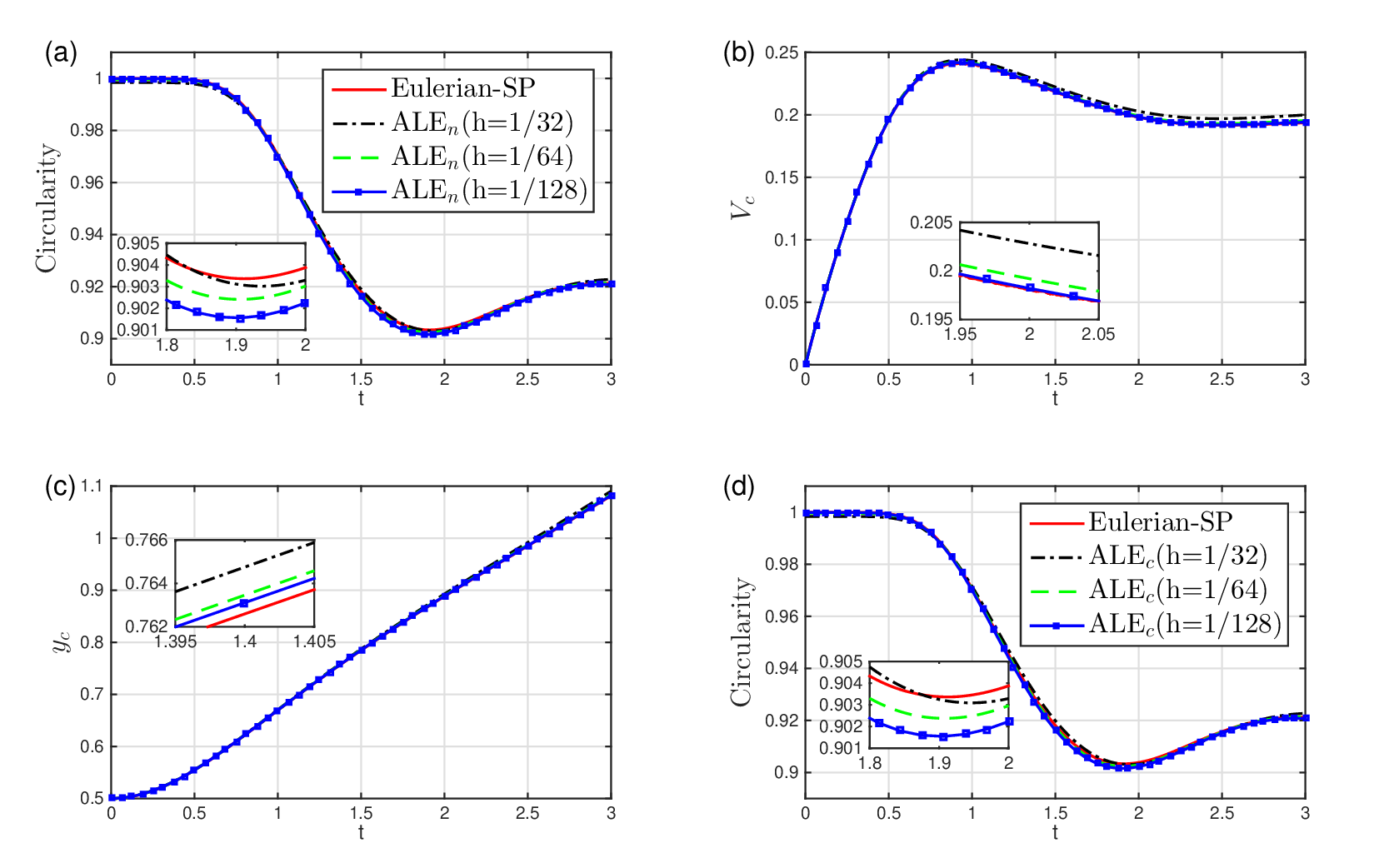}
\caption{The time history of the benchmark quantities for case I in 2d, where (a),(b),(c) are from the $\ALE_n$-SP method and (d) is from $\ALE_c$-SP method. The numerical results are compared with those from the \ESP-SP method using $2{\rm adapt}_{9,4}$.}
\label{fig:Fig3}
\end{figure}

\vspace{0.6em}\noindent
{\bf Example 4}: In this example, we again focus on case I and apply the two ALE structure-preserving methods to the rising bubble. We consider the P2-P0 element in \eqref{eq:P2P0} and the P2-(P1+P0) element in \eqref{eq:P2P1P0}, which both guarantee the assumption in \eqref{eq:ALEAsup2}. The benchmark results computed by the $\ALE_n$-SP and $\ALE_c$-SP methods are reported in Table \ref{tab:ALEBench}. Based on these observations, we can conclude that (i) the results from the two ALE methods are almost identical under the same computational parameters; (ii) the P2-(P1+P0) element yields more accurate and consistent results than the P2-P0 element. To further assess the performance of the two methods, we show the time evolution of the benchmark quantities in Fig.~\ref{fig:Fig3} and compare them with those from the \ESP-SP method using $2{\rm adapt}_{9,4}$. We observe the convergence of the two ALE methods as the mesh is refined.  Moreover, the ALE methods can produce results that are very consistent with those from the \ESP-SP method.

\begin{figure}[!htp]
\centering
\includegraphics[width=0.9\textwidth]
{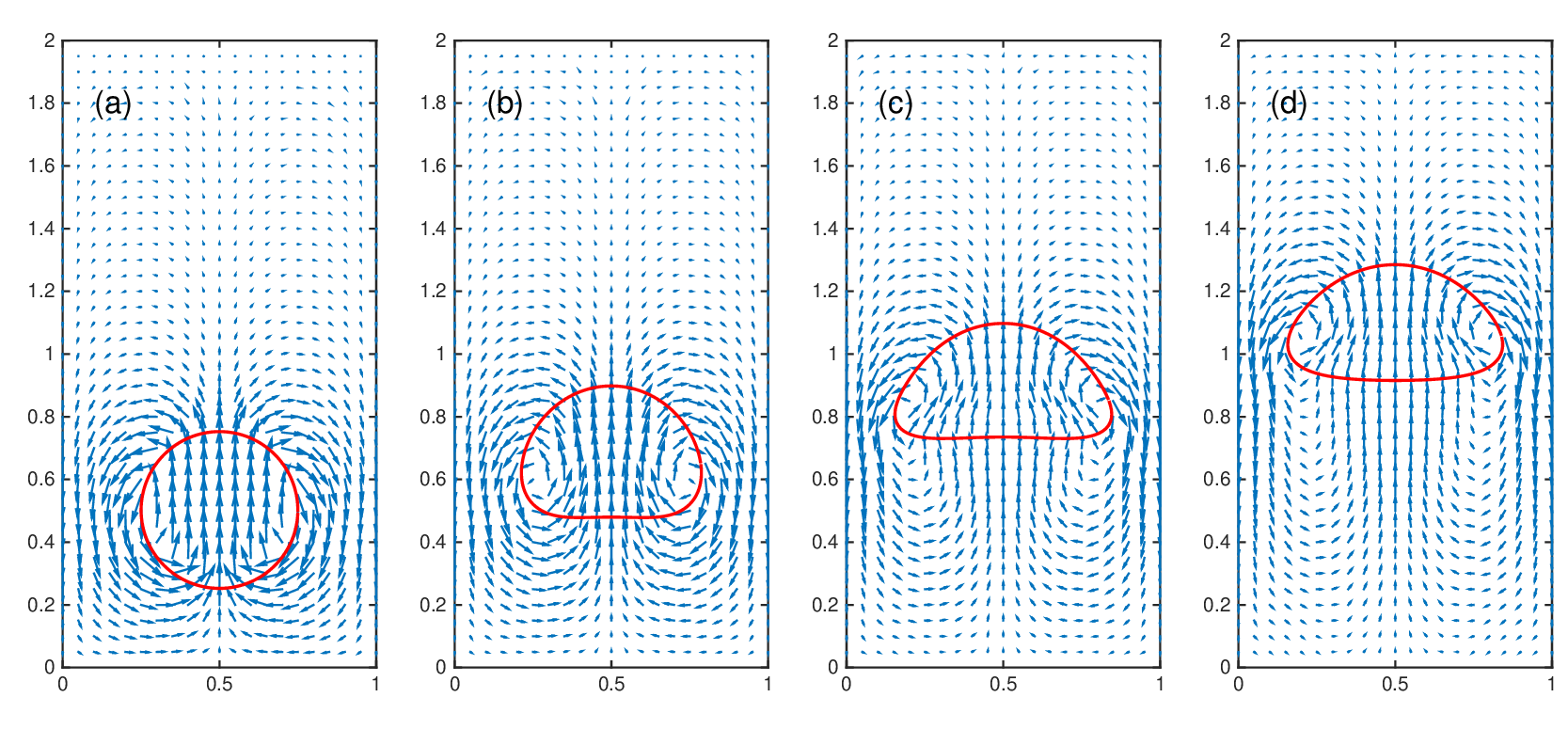}
\caption{Snapshots of the fluid interface (red line) and the velocity fields for the rising bubble in case I by using the $\ALE_n$-SP method, where $J_{\Gamma} = 128$, $\ttau=6.25\times 10^{-4}$, $J_\Omega=3794$. Here (a) $t=0.1$, $\max_{\vec x\in\Omega}|\vec u| = 0.0619$; (b) $t=1.0$, $\max_{\vec x\in\Omega}|\vec u| = 0.437$; (c) $t=2.0$, $\max_{\vec x\in\Omega}|\vec u| = 0.496$; (d) $t=3.0$; $\max_{\vec x\in\Omega}|\vec u| = 0.475$. }
\label{fig:4}
\end{figure}

\begin{figure}[!htp]
\centering
\includegraphics[width=0.9\textwidth]
{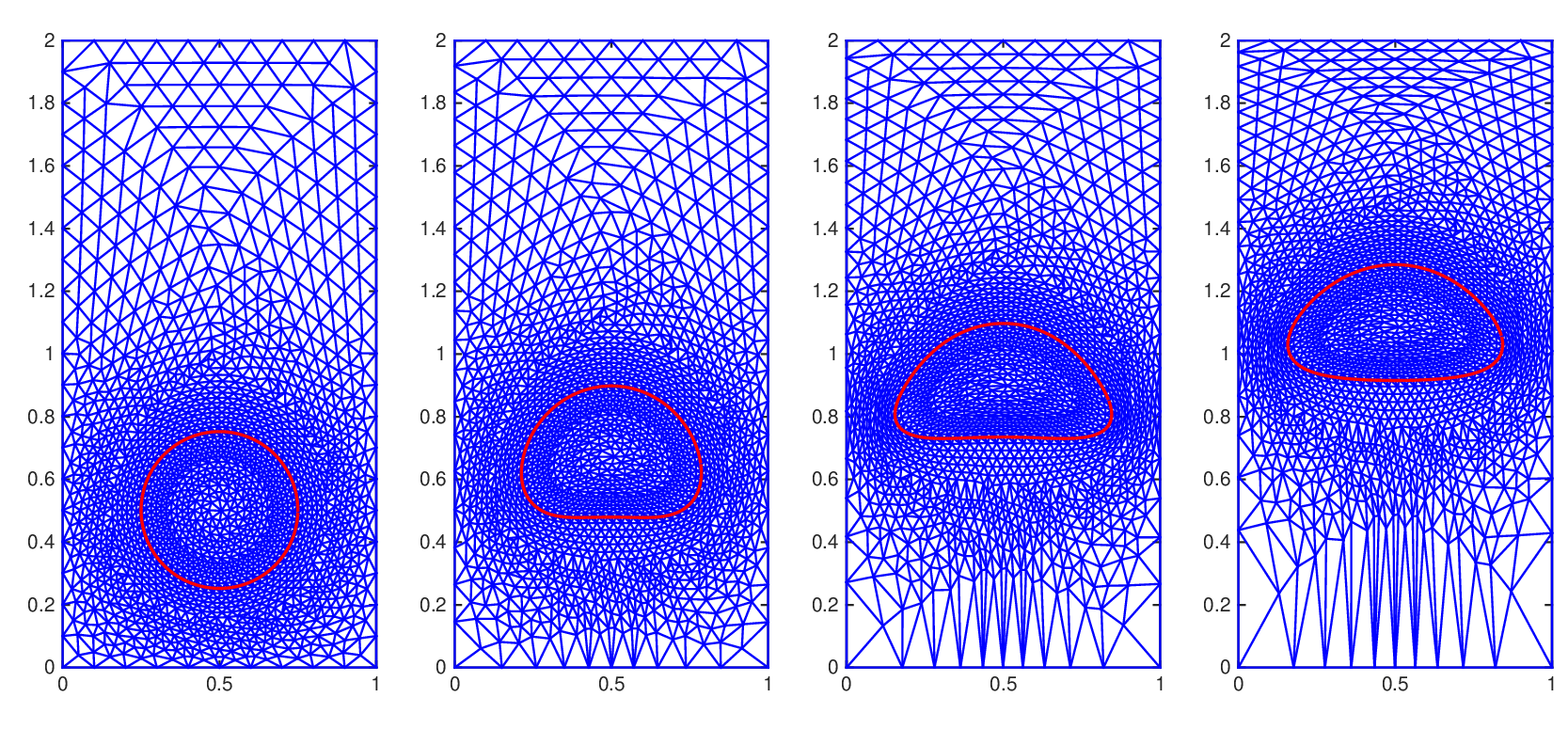}
\caption{ Snapshots of the corresponding computational meshes for the rising bubble in Fig.~\ref{fig:4}. }
\label{fig:5}
\end{figure}

\begin{figure}[!htp]
\centering
\includegraphics[width=0.8\textwidth]
{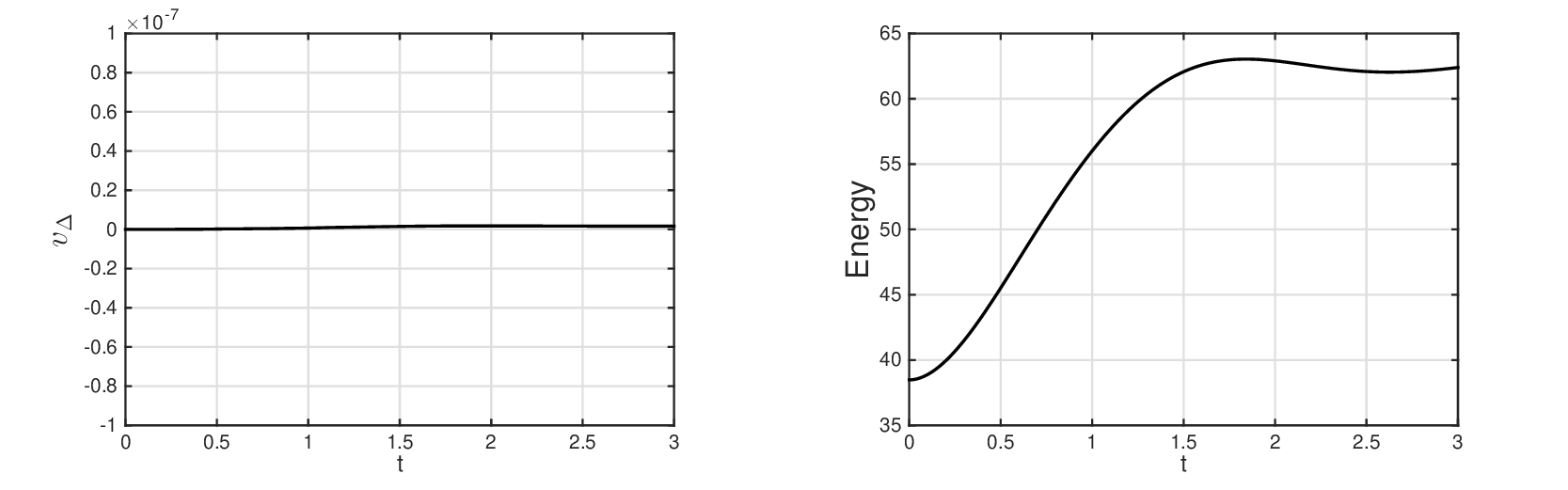}
\caption{Plots of the relative volume loss $v_\Delta$ and the discrete energy over time for the rising bubble in Fig.\ref{fig:4}. }
\label{fig:6}
\end{figure}

For the experiment using the $\ALE_n$-SP method with $J_\Gamma = 128$, $\ttau = 6.25\times 10^{-4}$, we show snapshots of the fluid interface and the velocity fields at several times in Fig.~\ref{fig:4}, and the corresponding computational meshes are presented in Fig.~\ref{fig:5}. We observe that the bulk mesh quality in the vicinity of the interface is generally well preserved. This shows that the moving mesh approach in section \eqref{sec:DALE} works quite smoothly, and no remeshing is necessary for this experiment.  The time history of the relative volume loss and the discrete energy are shown in Fig.~\ref{fig:6}. In particular, we observe the exact volume preservation.

\begin{figure}[!htp]
\centering
\includegraphics[width=0.8\textwidth]
{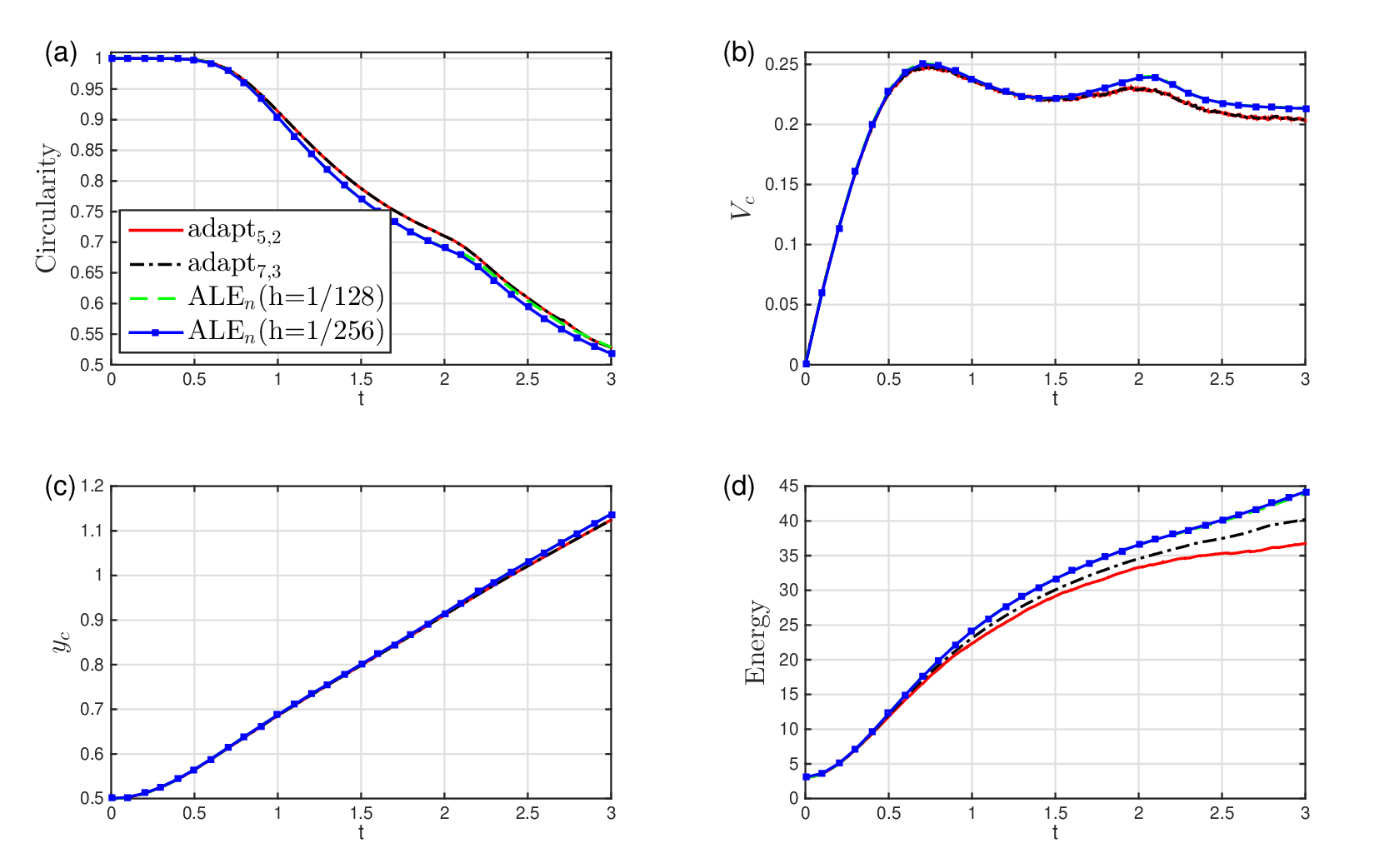}
\caption{Plots of (a) the circularity $\cira$, (b) the bubble's rise velocity $V_c$, (c) the bubble's center of mass $y_c$ and (d) the discrete energy for the rising bubble for case II in 2d by using the \ESP-SP and $\ALE_n$-SP methods.}\label{fig:6r}
\end{figure}

\vspace{0.6em}\noindent
{\bf Example 5}: In this example, we conduct a test of the rising bubble in case II, see \eqref{eq:Bench2}. The high ratio of the density and viscosity in the two phases will lead to a very strong deformation of the bubble. Therefore in the ALE moving mesh methods, the bulk mesh needs to be regenerated occasionally to facilitate the computation. In particular, we keep the interface mesh unchanged and regenerate the bulk mesh when the following condition is violated
\begin{equation}
\min_{\sigma\in\mathscr{T}^m}\min_{\alpha\in \measuredangle{(\sigma)}}\geq \tfrac1{18} \pi,
\end{equation}
where $\measuredangle{(\sigma)}$ is the set of all the angles of the simplex $\sigma$. The benchmark quantities for both the $\ESP$-SP and $\ALE_n$-SP methods are reported in Table \ref{tab:case2}. We also show the time plot of these benchmark values in Fig.~\ref{fig:6r}. Here we observe that the two methods can produce results that are very consistent with each other. 

\begin{figure}[!htp]
\centering
\includegraphics[width=0.85\textwidth]
{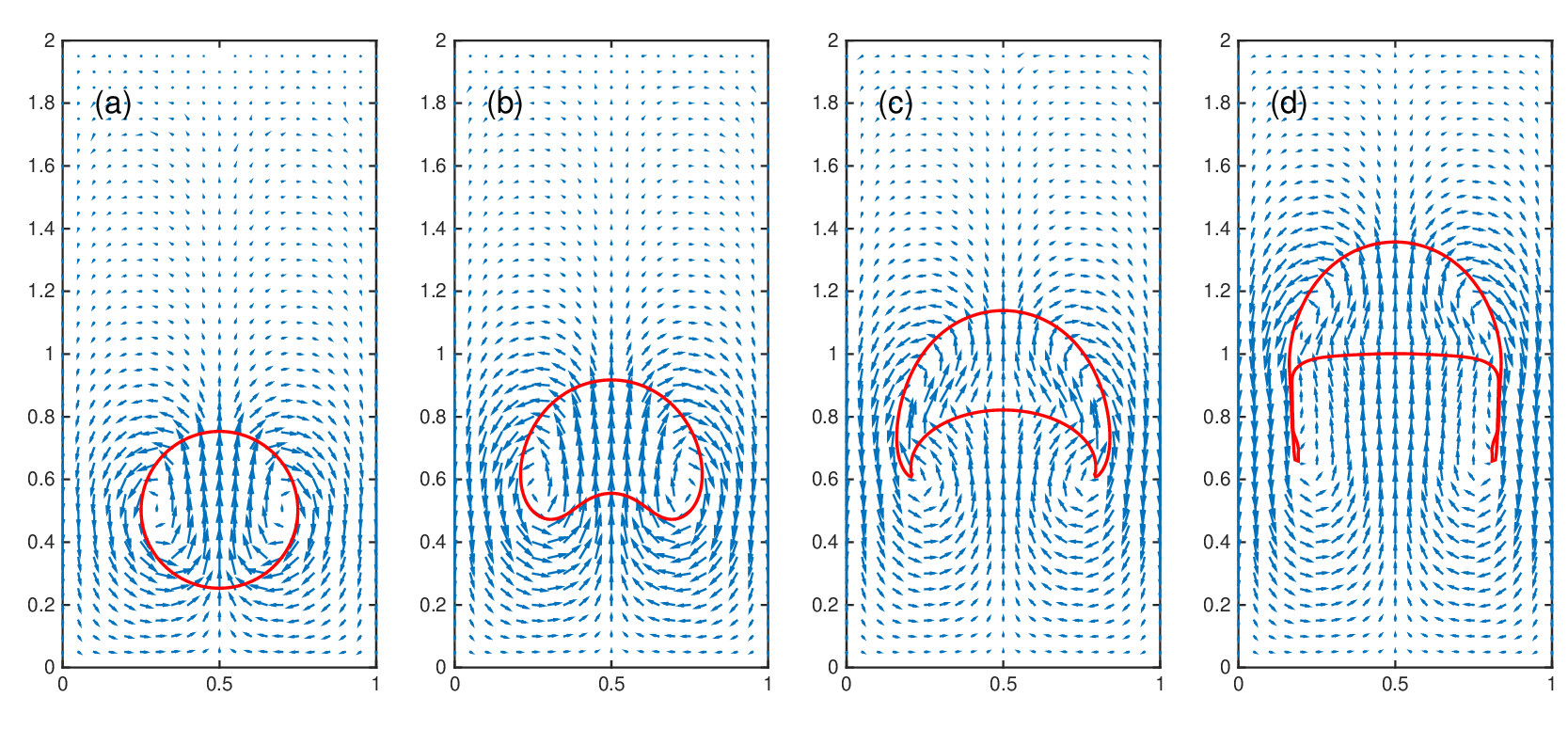}
\caption{ Snapshots of the fluid interface (red line) and the velocity fields for the rising bubble in case II by using the $\ALE_n$-SP method, where $J_{\Gamma} = 256$, $\ttau=1\times 10^{-3}$. Here (a) $t=0.1$, $\max_{\vec x\in\Omega}|\vec u| = 0.109$, $J_\Omega=13074$; (b) $t=1.0$, $\max_{\vec x\in\Omega}|\vec u| = 0.454$, $J_\Omega=14548$; (c) $t=2.0$, $\max_{\vec x\in\Omega}|\vec u| = 0.630$, $J_\Omega=17012$; (d) $t=3.0$; $\max_{\vec x\in\Omega}|\vec u| = 0.496$, $J_\Omega = 18688$.}\label{fig:7}
\end{figure}

\begin{figure}[!htp]
\centering
\includegraphics[width=0.8\textwidth]
{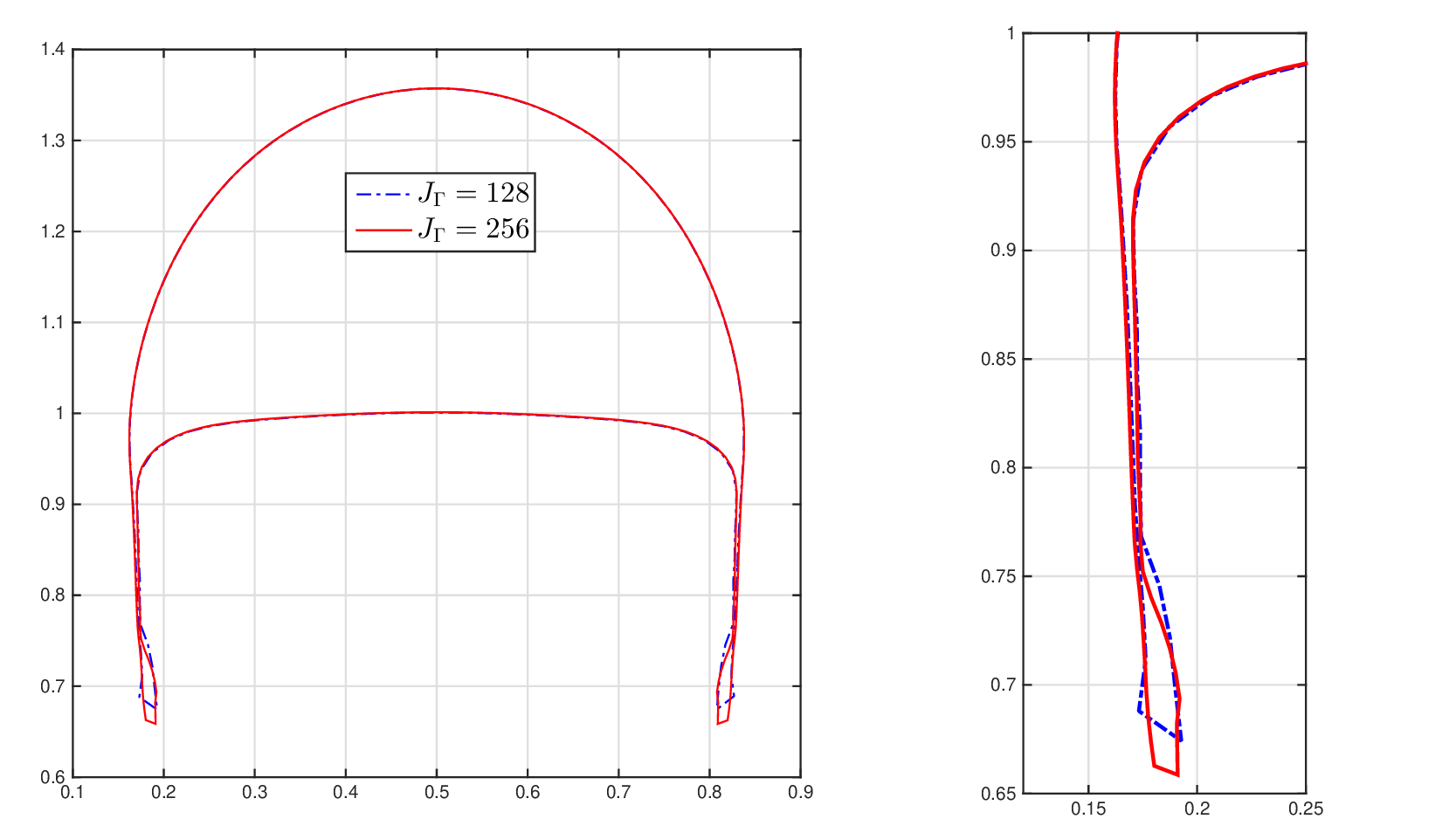}
\includegraphics[width=0.6\textwidth]
{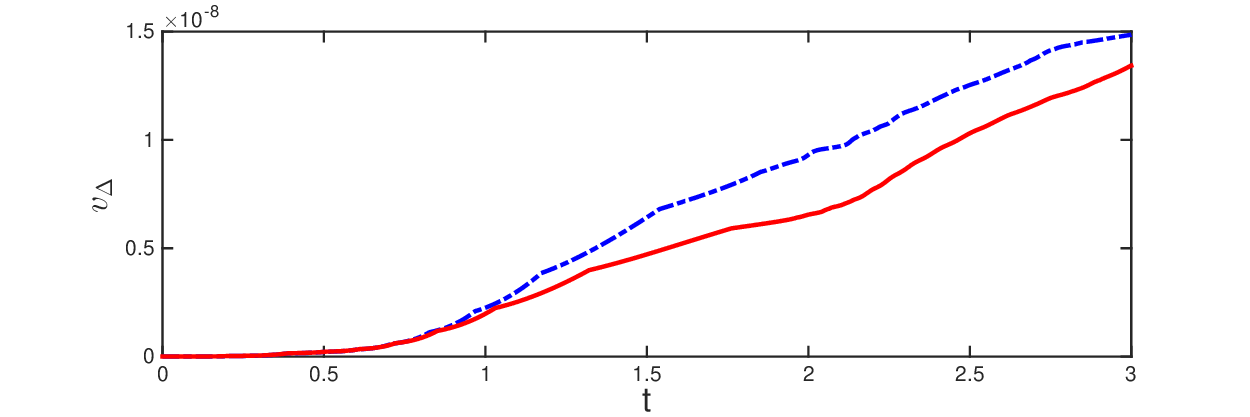}
\caption{Upper panel: Interface profiles at time $t=3$ by using the $\ALE_n$-SP method. Lower panel: the time history of the relative volume loss of the rising bubble.}
\label{fig:8}
\end{figure}

For the experiment using the $\ALE_n$-SP method with $J_\Gamma=256$, $\ttau = 10^{-3}$, we show snapshots of the interface and the velocity fields at several times  in Fig.~\ref{fig:7}. Moreover, the interface profiles at the final time $t=3$ for two different mesh sizes are depicted in Fig.~\ref{fig:8}. Here in both cases we observe very narrow tail structures in the bubble, and the interface profiles are quite close to each other. In particular, the volume preservation is observed as well in the lower panel of Fig.~\ref{fig:8}. 

\begin{table}[t]
\centering
\def\temptablewidth{0.7\textwidth}
\vspace{-2pt}
\caption{2d benchmark quantities of the rising bubble in case II. Here for the \ESP-SP method we use the pair element P2-P1 with XFEM, while for the $\ALE_n$-SP method we denote $h=1/J_\Gamma$ and use the pair element P2-(P1+P0) with $h_0=1/128, \ttau_0=10^{-3}$. Here ``${\rm Hysing3}$'' denotes the finest discretization run of group 3 in \cite{Hysing2009}.}\label{tab:case2}
{\rule{\temptablewidth}{1pt}}
\begin{tabular}{c|cc|cc|c}
 & \multicolumn{2}{c|}{Eulerian-SP} & \multicolumn{2}{c|}{${\rm ALE}_n$-SP} &${\rm Hysing3}$\\\hline 
  &${\rm adapt}_{5,2}$  & ${\rm adapt}_{7,3}$   &($h_0$, $\ttau_0)$ &($\frac{h_0}{2}$, $\ttau_0$) \\\hline  
$\cira_{\min}$ & 0.5885 &0.5282    &0.5285 &0.5180 &0.5144\\\hline 
$t_{\cira = \cira_{\min}}$ &3.0000 &3.0000  &3.0000 & 3.0000 &3.0000\\[0.4em]\hline 
$V_{c,\max 1}$ &0.2583 &0.2480  & 0.2503 & 0.2502 &0.2502\\\hline 
$t_{_{V_c = V_{c,\max 1}}}$ &0.8800 &0.7600    &0.7300 & 0.7300 &0.7317\\\hline 
$V_{c,\max 2}$ &0.2286 &0.2306   &0.2400 &0.2398 &0.2393\\\hline 
$t_{_{V_c = V_{c,\max 2}}}$ &2.0000 &1.9520   &2.0690 & 2.0600 &2.0600\\\hline 
$y_c(t=3)$ &1.1274 &  1.1243   &1.1379 &1.1377 &1.1376
\end{tabular}
{\rule{\temptablewidth}{1pt}}
\end{table}

\subsubsection{Numerical results in 3d}

\begin{table}
\centering
\def\temptablewidth{0.70\textwidth}
\vspace{-4pt}
\caption{3d benchmark quantities of the rising bubble for structure-preserving methods in case I. Here for the \ESP-SP method we use the pair element P2-P1 with XFEM, while for the $\ALE_n$-SP method we use the pair element P2-(P1+P0) and $J_\Gamma=842$, $\ttau_0=5\times 10^{-3}$. Here ``BGN'' denotes the finest discretization run from \cite{BGN15stable}.}
\label{tab:3db1}
{\rule{\temptablewidth}{1pt}}
\begin{tabular}{c|cc|cc|c}
 & \multicolumn{2}{c|}{Eulerian-SP} & \multicolumn{2}{c|}{${\rm ALE}_n$-SP} &BGN\\\hline 
   &${\rm adapt}_{5,2}$  &${\rm adapt}_{6,3}$   &$\ttau=\ttau_0$ &$\ttau=\frac{\ttau_0}{5}$ &${\rm adapt}_{6,3}$ \\\hline  
$\cira_{\min}$  & 0.9570  &0.9507   &0.9422 &0.9422 &0.9508\\\hline 
$t_{\cira = \cira_{\min}}$  &3.0000 &3.0000 &3.0000 & 3.0000 &3.0000\\[0.4em]\hline 
$V_{c,\max}$  &0.3821  &0.3846   & 0.3949 & 0.3945 &0.3845\\\hline 
$t_{_{V_c = V_{c,\max}}}$  &1.1690 &1.0800   &0.9950 & 0.9980 &1.0800\\\hline 
$y_c(t=3)$  &1.5516 &1.5555  &1.5710 &1.5702 &1.5555
\end{tabular}
{\rule{\temptablewidth}{1pt}}
\end{table}

\begin{figure}[!htp]
\centering
\includegraphics[width=0.8\textwidth]{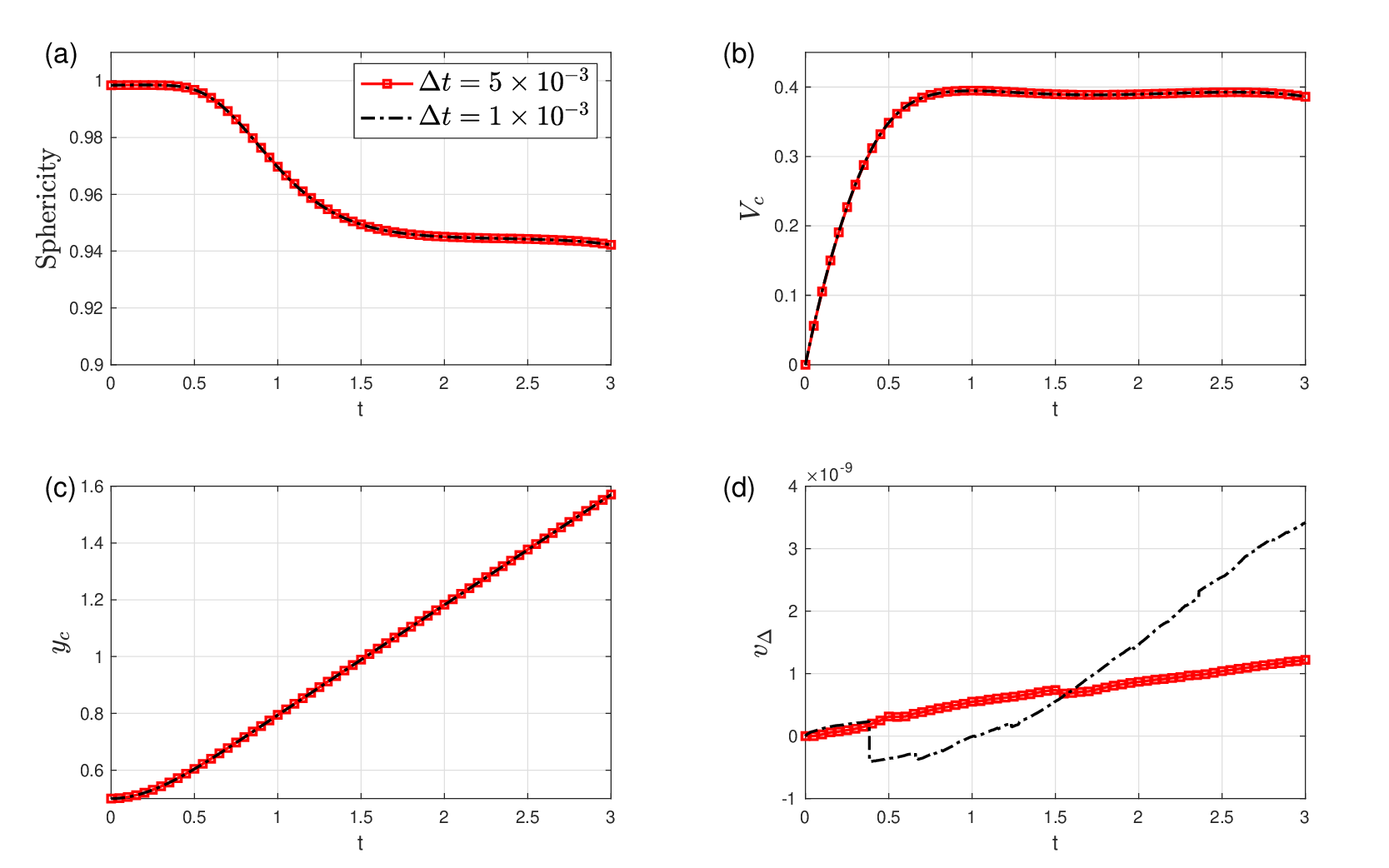}
\caption{The time history of the sphericity, the rise velocity, the center of mass and the relative volume loss of the rising bubble for case I in 3d, here the results are obtained from the $\ALE_n$-SP method with $J_\Gamma=842$.}
\label{fig:10} 
\end{figure}

\begin{figure}[!htp]
\centering
\includegraphics[width=0.8\textwidth]{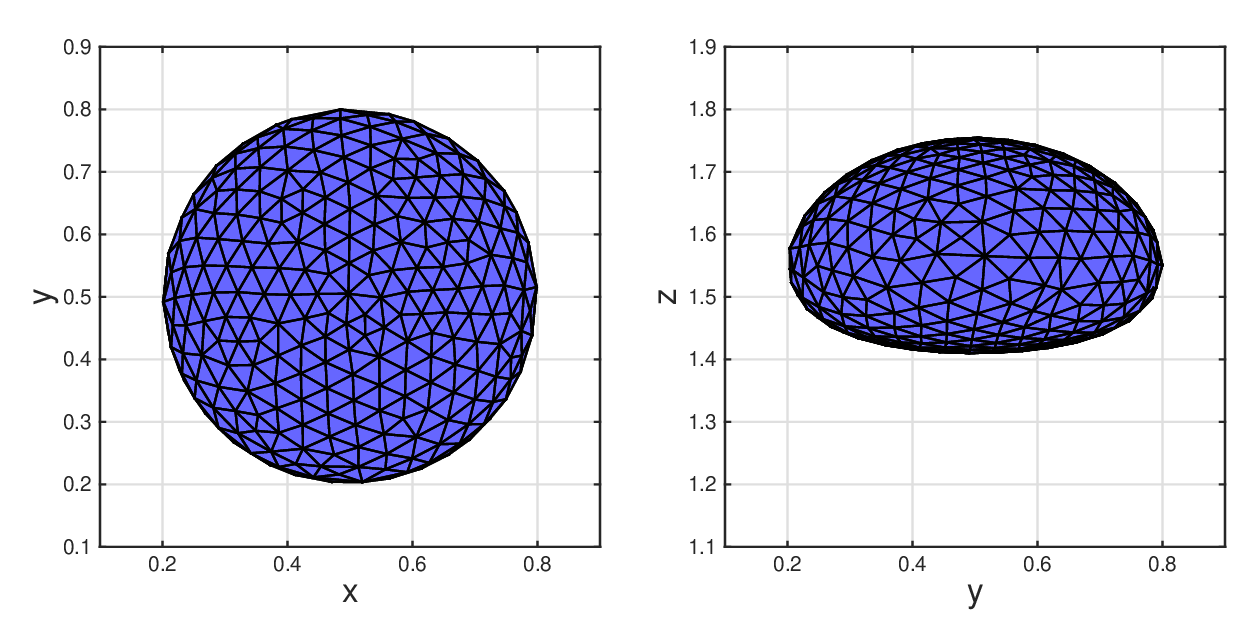}
\caption{The interface mesh for the 3d rising bubble of case I at time $T = 3$. Left panel: view from the top; right panel: view from the front.}
\label{fig:11}
\end{figure}

\vspace{0.4em}\noindent
{\bf Example 6}: 
We consider the natural three-dimensional analogue of the
computational setup from \S\ref{sec:nr2d}. That is, we let
 $\overline{\Omega} = [0,1]\times[0,1]\times[0,2]$ with $\partial_1\Omega=[0,1]\times[0,1]\times\{0,2\}$, and $\partial_2\Omega=\partial\Omega\backslash\partial_1\Omega$. The initial interface is given by $\Gamma(0)=\bigl\{\vecz\in\Omega:\;|\vecz - (\frac{1}{2}, \frac{1}{2}, \frac{1}{2})^T|=\frac{1}{4}\bigr\}$.
In this example, we are focused on the rising bubble in case I.
The quantitative values for the rising bubble are reported in Table \ref{tab:3db1} for both the \ESP-SP and $\ALE_n$-SP methods, where we observe a good agreement. For the experiments using the $\ALE_n$-SP method, we also plot the time history of the benchmark quantities in Fig.~\ref{fig:10}. We observe the volume preservation in Fig.~\ref{fig:10}(d). Moreover, visualizations of the interface mesh for the final bubble are shown in Fig.~\ref{fig:11}. Here the unfitted results are near to that in \cite{BGN15stable}, where the sphericity decreases from 1 to around 0.95. For the ALE method, the results are quite similar to that in \cite{Duan2022energy}. Moreover, we expect that the \ESP-SP method can produce more closer results to those of the $\ALE_n$-SP method if a finer mesh size is employed.

\vspace{0.4em}\noindent
{\bf Example 7}: 
In order to be able to perform a quantitative comparison with the computations in \cite{adelsberger20143d}, we also briefly consider the exact same setup from
Example 6, but now with no-slip boundary conditions on all of $\partial\Omega$,
i.e.\ $\partial_1\Omega=\partial\Omega$. The obtained benchmark results are reported in Table~\ref{tab:3db2}, where the presented quantities from \cite{adelsberger20143d} were extracted visually from their graphs to serve as a comparison. We observe these results show a good agreement, which further verifies the accuracy of our numerical methods. 


\begin{table}[!htp]
\centering
\def\temptablewidth{0.75\textwidth}
\vspace{-4pt}
\caption{3d benchmark quantities of the rising bubble in case I for structure-preserving methods with the same setting in Table \ref{tab:3db1} except that no-slip boundary conditions are employed on all of $\partial\Omega$.  Here ``Adelsberger'' denotes the discretization run from \cite{adelsberger20143d}.}
\label{tab:3db2}
{\rule{\temptablewidth}{1pt}}
\begin{tabular}{c|cc|ccc}
 & Eulerian-SP & ${\rm ALE}_n$-SP &\multicolumn{3}{c}{{\rm Adelsberger}}\\\hline 
   &${\rm adapt}_{6,3}$  &$\ttau=\ttau_0$ &Drops & NaSt3DGPF & OpenFoam \\\hline  
$\cira_{\min}$  &  0.9657  &0.9605   &0.9600 &0.9600 &0.9550\\\hline 
$t_{\cira = \cira_{\min}}$  & 1.9770 &1.7950 &2.1500 & 1.8500 &2.0000\\[0.4em]\hline 
$V_{c,\max}$  &0.3519  & 0.3603   & 0.3570 & 0.3585 &0.3520\\\hline 
$t_{_{V_c = V_{c,\max}}}$  & 0.7760 &0.8800   &0.8400 & 0.8400 &0.8400\\\hline 
$y_c(t=3)$  & 1.4607 &1.4680  &1.4750 &1.4700 &1.4300
\end{tabular}
{\rule{\temptablewidth}{1pt}}
\end{table}

\section{Conclusion}\label{sec:con}
In this work, we proposed two structure-preserving methods for discretizing 
two-phase Navier--Stokes flow using either an unfitted or a fitted mesh 
approach. 
The proposed methods combine a parametric finite element approximation for the 
evolving interface together with unfitted or fitted finite element 
approximations for the Navier--Stokes equations in the bulk. The parametric 
approximation is based on the BGN formulation, which allows for tangential 
degrees of freedoms and leads to a good mesh quality of the interface 
approximation. The unfitted and fitted approximations are based on an 
Eulerian and ALE weak formulation, respectively. We proved that the 
resulting two methods satisfy unconditional stability and exact volume 
preservation. Numerical results were presented to demonstrate the accuracy and 
efficiency of the introduced methods, and to numerically verify these 
structure-preserving properties.

All in all, the numerical results from the unfitted and fitted mesh approaches are quite similar. In several numerical examples, the latter approach seems to be able to provide more accurate results, in particular in terms of the accuracy of the pressure approximations. This is due to the less accurate capture of the jumps across the interface in the physical parameters and in the solutions. Nevertheless, the unfitted mesh approach is still desirable especially in the case when the interface may undergo large deformations and topological changes. Further investigations of the introduced structure-preserving unfitted method in the context of cutFEM is left for future research.

\section*{Acknowledgement}
The work of Quan Zhao was funded by the Alexander von Humboldt Foundation. 

\appendix

\section{Derivation of \eqref{eq:cALEweak1}}\label{app:nc}

Let $\varphi: \Omega\times[0,T]\to\bR$ be a scalar field. In terms of the ALE reference domain, applying the Reynolds transport theorem yields that
\begin{equation}
\ddt\int_{\Omega_\pm(t)}\varphi\,\dL^d  = \int_{\Omega_\pm(t)}(\partial_t^\circ\varphi + \varphi\,\nabla\cdot\vec w)\,\dL^d,\label{eq:ALERey}
\end{equation}
where $\partial_t^\circ$ is the derivative with respect to the ALE reference, see \eqref{eq:ALED}, and $\vec w$ is the corresponding ALE frame velocity given in \eqref{eq:ALEmeshv}.  

Combining \eqref{eq:ALERey} with the product rule for $\partial_t^\circ$
yields that
\begin{equation}
\ddt\bigl(\vec u,~\vec\chi\bigr)_{\Omega_\pm(t)} = \bigl(\partial_t^\circ\vec u,~\vec\chi\bigr)_{\Omega_\pm(t)} + \bigl(\vec u,~\partial_t^\circ\vec\chi\bigr)_{\Omega_\pm(t)} + \bigl(\vec u\cdot\vec\chi,~\nabla\cdot\vec w\bigr)_{\Omega_\pm(t)}\quad\forall\vec\chi\in\mathbb{V},
\label{eq:ALEuv}
\end{equation}
where $\mathbb{V}$ is the function space defined in \eqref{eqn:UPspaces}. Multiplying \eqref{eq:ALEuv} with $\rho_\pm$ and summing gives
\begin{equation}
\ddt\bigl(\rho\,\vec u,~\vec\chi\bigr)-\tfrac{1}{2}\bigl(\rho\,\nabla\cdot\vec w,~\vec u\cdot\vec\chi\bigr) - \bigl(\rho\,\vec u,~\partial_t^\circ\vec\chi\bigr)=\bigl(\rho\,\partial_t^\circ\vec u,~\vec\chi\bigr)+\tfrac{1}{2}\bigl(\rho\,\nabla\cdot\vec w,~\vec u\cdot\vec\chi\bigr),
\end{equation}
which then implies \eqref{eq:cALEweak1} by comparing with \eqref{eq:ALEweak1}.

\bibliographystyle{model1b-num-names}
\bibliography{bib}
\end{document}